
\input amssym.def
\input amssym.tex


\def\item#1{\vskip1.3pt\hang\textindent {\rm #1}}


\tolerance=300
\pretolerance=200
\hfuzz=1pt
\vfuzz=1pt


\hoffset=0.6in
\voffset=0.8in

\hsize=5.8 true in 
\vsize=9.2 true in
\parindent=25pt
\mathsurround=1pt
\parskip=1pt plus .25pt minus .25pt
\normallineskiplimit=.99pt

\countdef\revised=100
\mathchardef\emptyset="001F 
\chardef\ss="19
\def\3{\ss}
\def\anf{$\lower1.2ex\hbox{"}$}
\def\frac#1#2{{#1 \over #2}}
\def\>{>\!\!>}
\def\<{<\!\!<}

\def\ssarr{\hbox to 30pt{\rightarrowfill}}
\def\sarr{\hbox to 40pt{\rightarrowfill}}
\def\arr{\hbox to 60pt{\rightarrowfill}}
\def\larr{\hbox to 60pt{\leftarrowfill}}
\def\Arr{\hbox to 80pt{\rightarrowfill}}

{}

\def\ad{\mathop{\rm ad}\nolimits}

\def\Ad{\mathop{\rm Ad}\nolimits}

\def\Aut{\mathop{\rm Aut}\nolimits}

\def\cone{\mathop{\rm cone}\nolimits}
\def\conv{\mathop{\rm conv}\nolimits}
\def\Der{\mathop{\rm Der}\nolimits}

\def\Exp{\mathop{\rm Exp}\nolimits}

\def\End{\mathop{\rm End}\nolimits}
\def\Ext{\mathop{\rm Ext}\nolimits}

\def\Hol{\mathop{\rm Hol}\nolimits}%
\def\id{\mathop{\rm id}\nolimits} 

\def\inf{\mathop{\rm inf}\nolimits}

\def\Inn{\mathop{\rm Inn}\nolimits}
\def\Int{\mathop{\rm int}\nolimits}

\def\Prim{\mathop{\rm Prim}\nolimits}



\def\Sl{\mathop{\rm Sl}\nolimits}

\def\span{\mathop{\rm span}\nolimits}

\def\Spec{\mathop{\rm Spec}\nolimits}

\def\sup{\mathop{\rm sup}\nolimits}

\def\tr{\mathop{\rm tr}\nolimits}

\def\0{{\bf 0}}
\def\1{{\bf 1}}

\def\a{{\frak a}}

\def\g{{\frak g}}

\def\k{{\frak k}}

\def\m{{\frak m}}

\def\n{{\frak n}}

\def\p{{\frak p}}

\def\s{{\frak s}}

\def\sL{{\frak {sl}}}
\def\t{{\frak t}}

\def\z{{\frak z}}

\def\C{{\Bbb C}}

\def\N{{\Bbb N}}

\def\R{{\Bbb R}}

\def\:{\colon}  
\def\.{{\cdot}}
\def\|{\Vert}
\def\bsk{\bigskip}

\def\giantskip{\vskip2\bigskipamount}
\def\gsk{\giantskip}
\def \la {\langle}
\def\msk{\medskip}
\def \ra {\rangle}
\def \res {\!\mid\!\!}

\def\bbr{\bigbreak}
\def\giantbreak{\par \ifdim\lastskip<2\bigskipamount \removelastskip
         \penalty-400 \giantskip\fi}

\def\nin{\noindent}
\def\cen{\centerline}
\def\pagebreak{\vskip 0pt plus 0.0001fil\break}
\def\linebreak{\break}

\def\hat{\widehat}

\def\eps{\varepsilon}
\def\epsilon{\varepsilon}
\def\eset{\emptyset}

\def\nin{\noindent}
\def\oline{\overline}

\def\pder#1,#2,#3 { {\partial #1 \over \partial #2}(#3)}
\def\pde#1,#2 { {\partial #1 \over \partial #2}}
\def\phi{\varphi}


\def\subeq{\subseteq}
\def\supeq{\supseteq}

\def\tilde{\widetilde}

\font\eightrm=cmr8


\font\bfone=cmbx10 scaled\magstep1 
\font\bftwo=cmbx10 scaled\magstep2 

\def\qed{{\unskip\nobreak\hfil\penalty50\hskip .001pt \hbox{}\nobreak\hfil
          \vrule height 1.2ex width 1.1ex depth -.1ex
           \parfillskip=0pt\finalhyphendemerits=0\medbreak}\rm}


\def\Lemma #1. {\bigbreak\vskip-\parskip\noindent{\bf Lemma #1.}\quad\it}

\def\Sublemma #1. {\bigbreak\vskip-\parskip\noindent{\bf Sublemma #1.}\quad\it}

\def\Proposition #1. {\bigbreak\vskip-\parskip\noindent{\bf Proposition #1.}
\quad\it}

\def\Corollary #1. {\bigbreak\vskip-\parskip\nin{\bf Corollary #1.}
\quad\it}

\def\Theorem #1. {\bigbreak\vskip-\parskip\noindent{\bf Theorem #1.}
\quad\it}

\def\Definition #1. {\rm\bigbreak\vskip-\parskip\noindent{\bf Definition #1.}
\quad}

\def\Remark #1. {\rm\bigbreak\vskip-\parskip\noindent{\bf Remark #1.}\quad}

\def\Example #1. {\rm\bigbreak\vskip-\parskip\noindent{\bf Example #1.}\quad}

\def\Problems #1. {\bigbreak\vskip-\parskip\noindent{\bf Problems #1.}\quad}
\def\Problem #1. {\bigbreak\vskip-\parskip\noindent{\bf Problems #1.}\quad}

\def\Conjecture #1. {\bigbreak\vskip-\parskip\noindent{\bf Conjecture #1.}\quad}

\def\Proof#1.{\rm\par\ifdim\lastskip<\bigskipamount\removelastskip\fi\smallskip
            \noindent {\bf Proof.}\quad}

\def\Axiom #1. {\bigbreak\vskip-\parskip\noindent{\bf Axiom #1.}\quad\it}

\def\Satz #1. {\bigbreak\vskip-\parskip\noindent{\bf Satz #1.}\quad\it}

\def\Korollar #1. {\bbr\vskip-\parskip\nin{\bf Korollar #1.} \quad\it}

\def\Bemerkung #1. {\rm\bigbreak\vskip-\parskip\noindent{\bf Bemerkung #1.}
\quad}

\def\Beispiel #1. {\rm\bigbreak\vskip-\parskip\noindent{\bf Beispiel #1.}\quad}
\def\Aufgabe #1. {\rm\bigbreak\vskip-\parskip\noindent{\bf Aufgabe #1.}\quad}

\def\Beweis#1. {\rm\par\ifdim\lastskip<\bigskipamount\removelastskip\fi
           \smallskip\noindent {\bf Beweis.}\quad}

\nopagenumbers

\def\date{\ifcase\month\or January\or February \or March\or April\or May
\or June\or July\or August\or September\or October\or November
\or December\fi\space\number\day, \number\year}

\def\title{Title ??}
\def\author{Author ??}

\def\thanks#1{\footnote*{\eightrm#1}}

\def\rightheadline{\hfil{\eightrm\title}\hfil\tenbf\folio}
\def\leftheadline{\tenbf\folio\hfil{\eightrm\author}\hfil}
\headline={\vbox{\line{\ifodd\pageno\rightheadline\else\leftheadline\fi}}}

\def\firstheadline{}
\def\firstfootline{\cen{\rm\folio}}

\def\seite #1 {\pageno #1
               \headline={\ifnum\pageno=#1 \firstheadline
               \else\ifodd\pageno\rightheadline\else\leftheadline\fi\fi}
               \footline={\ifnum\pageno=#1 \firstfootline\else{}\fi}}

\newdimen\dimenone
 \def\checkleftspace#1#2#3#4{
 \dimenone=\pagetotal
 \advance\dimenone by -\pageshrink   
 \ifdim\dimenone>\pagegoal          
   \else\dimenone=\pagetotal
        \advance\dimenone by \pagestretch
        \ifdim\dimenone<\pagegoal
          \dimenone=\pagetotal
          \advance\dimenone by#1         
          \setbox0=\vbox{#2\parskip=0pt                
                     \hyphenpenalty=10000
                     \rightskip=0pt plus 5em
                     \noindent#3 \vskip#4}    
        \advance\dimenone by\ht0
        \advance\dimenone by 3\baselineskip   
        \ifdim\dimenone>\pagegoal\vfill\eject\fi
          \else\eject\fi\fi}


\def\subheadline #1{\nin\bigbreak\vskip-\lastskip
      \checkleftspace{0.7cm}{\bf}{#1}{\medskipamount}
          \indent\vskip0.7cm\centerline{\bf #1}\medskip}

\def\sectionheadline #1{\bigbreak\vskip-\lastskip
      \checkleftspace{1.1cm}{\bf}{#1}{\bigskipamount}
         \vbox{\vskip1.1cm}\cen{\bfone #1}\bsk}

\def\lsectionheadline #1 #2{\bigbreak\vskip-\lastskip
      \checkleftspace{1.1cm}{\bf}{#1}{\bigskipamount}
         \vbox{\vskip1.1cm}\cen{\bfone #1}\msk \cen{\bfone #2}\bsk}

\def\lchapterheadline #1 #2{\bigbreak\vskip-\lastskip\indent\vskip3cm
                       \cen{\bftwo #1} \msk \cen{\bftwo #2} \gsk}
\def\llsectionheadline #1 #2 #3{\bigbreak\vskip-\lastskip\indent\vskip1.8cm
\cen{\bfone #1} \msk \cen{\bfone #2} \msk \cen{\bfone #3} \nobreak\bsk\nobreak}


\newtoks\literat
\def\[#1 #2\par{\literat={#2\unskip.}%
\hbox{\vtop{\hsize=.15\hsize\nin [#1]\hfill}
\vtop{\hsize=.82\hsize\nin\the\literat}}\par
\vskip.3\baselineskip}

\mathchardef\emptyset="001F 
\def\address{Author: \tt$\backslash$def$\backslash$address$\{$??$\}$}

\def\firstpage{\nin
{\obeylines \parindent 0pt }
\vskip2cm
\centerline {\bfone \title}
\gsk
\centerline{\bf\author}

\vskip1.5cm \rm}

\def\addresstwo{}

\def\dlastpage{\par\vbox{\vskip1cm\nin
\line{
\vtop{\hsize=.5\hsize{\parindent=0pt\baselineskip=10pt\nin\address}}
\quad 
\vtop{\hsize=.42\hsize\nin{\parindent=0pt
\baselineskip=10pt\addresstwo}}
\hfill} }}


\magnification=\magstep0

\pageno=1

\def\title{On the dual of complex Ol'shanski\u\i{} semigroups}
\def\author{Bernhard Kr\"otz${}^*$}
\footnote{}{${}^*$Supported by the DFG-project HI 412/5-2}
\def\date{November 10, 1999}
\def\leftheadline{\tenbf\folio\hfil\eightrm\date}
\def\bs{\backslash} 
\def\Prim{\mathop{\rm Prim}\nolimits}
\def\address
{Bernhard Kr\"otz

The Ohio State University

Department of Mathematics 

231 West 18th Avenue

Columbus, OH 43210-1174

USA}

\firstpage

\sectionheadline{Introduction}

Let $G$ be a connected Lie group which sits in its 
universal complexification $G_\C$.  If $\g$ denotes the Lie algebra
of $G$ and $W\subeq \g$ is a
non-empty open $\Ad(G)$-invariant convex cone
containing no affine lines, then a {\it complex
Ol'shanski\u\i{} semigroup} is defined by $S=G\exp(iW)\subeq G_\C$. 
One knows that $S$ is an open subsemigroup of $G_\C$ with holomorphic
multiplication, and moreover $S$ is invariant under the
antiholomorphic involution $g\mapsto g^*=\oline g^{-1}$, where $\oline
g$ indicates complex conjugation in $G_\C$. 
In particular, $(S,*)$ is an {\it involutive semigroup}. 

\par If ${\cal H}$ is a Hilbert space and $B({\cal H})$ denotes the bounded
operators on it, then a {\it holomorphic representation}
$(\pi, {\cal H})$ of $S$ is a holomorphic semigroup homomorphism 
$\pi\: S\to B({\cal H})$ with total image and which satisfies 
$\pi(s^*)=\pi(s)^*$ for all $s\in S$.  
A mapping $\alpha\: S\to [0,\infty[$ satisfying $\alpha(s^*)=\alpha(s)$ and
$\alpha(st)\leq \alpha(s)\alpha(t)$  for all $s,t\in S$ is called an
{\it absolute value} of $S$. We call a  holomorphic representation
$\alpha$-{\it bounded} for some absolute value $\alpha$ 
if $\|\pi(s)\|\leq \alpha(s)$ holds for $s\in S$.

\par The $\alpha$-bounded holomorphic representations of $S$ can be
modelled via a certain $C^*$-algebra  $C_h^*(S,\alpha)$  (cf.\
Definition I.3, Lemma II.6), i.e., there is a natural correspondence between
$\alpha$-bounded holomorphic representations of $S$ and non-degenerate
representations of $C_h^*(S,\alpha)$. An important result of K.-H. Neeb 
asserts that these $C^*$-algebras $C_h^*(S,\alpha)$ are CCR (cf.\
[Ne99, Ch.\ XI]). 
If we denote by $\hat S_\alpha$ the set of equivalence classes of 
irreducible $\alpha$-bounded representations of $S$, we therefore
obtain a bijection $\hat S_\alpha\cong C_h^*(S,\alpha)\, \hat{}$. The
hull-kernel topology on $C^*(S,\alpha)\,\hat{}$ thus gives rise to a topology
on $\hat S_\alpha$ denoted by ${\cal T}_{hk}^\alpha$.

\par On the other hand irreducible holomorphic representations of $S$
are obtained in a unique fashion by  analytic continuation of unitary
highest weight representations of $G$ and vice versa. Thus we may
think of $\hat S_\alpha$ as a certain subset of $\hat G$ and the
topology on $\hat G$ induces a topology on $\hat S_\alpha$ denoted by 
${\cal T}_G^\alpha$. Finally the parametrization of $\hat S_\alpha$ by
a certain subset of highest weights $HW_\alpha\subeq i\t^*$, where $\t$ denotes a compactly
embedded Cartan subalgebra, gives a topology on ${\cal T}_e^\alpha$
obtained by the euclidean topology on $HW_\alpha$.  Our main result
(cf.\ Theorem II.24) then asserts that one 
has  
$${\cal T}_{hk}^\alpha\subeq {\cal T}_G^\alpha\subeq{\cal T}_e^\alpha.$$   
Moreover, the induced Borel structures are all the same and for 
generic absolute values even all the topologies coincide. 
These  results are obtained by a combination of holomorphic representation
theory (cf.\ [Ne99]), standard structure theory of $C^*$-algebras 
(cf.\ [Dix82]) and real analysis methods 
(boundary values of Poisson transforms). 
Our results
imply in particular that the CCR-algebra $C_h^*(S,\alpha)$ has separated dual
for generic absolute values, which allows us to identify 
$C_h^*(S,\alpha)$ with the $C^*$-algebra 
defined by the continuous field of elementary $C^*$-algebras
$\big({\cal K}({\cal H}_\lambda)\big)_{\lambda\in HW_\alpha}$. Further we
explain what our result means for the 
abstract representation theory of complex Ol'shanski\u\i{}
semigroups. 

\par Finally we give a criterion for an absolute value $\alpha$ in order 
that the $C^*$-algebra $C_h^*(S,\alpha)$ has continuous trace (cf.\
Proposition II.33).

\sectionheadline{Involutive semigroups}

\Definition I.1. (Involutive semigroups) (a) An {\it involutive
semigroup} is a semigroup $S$
equipped with an involutive antiautomorphism $s\mapsto s^*$.
\par If $S$ does not have an identity, we set
$S_\1\:=S\dot\cup\{\1\}$. Then the prescription $s\1=\1s=s$ for all
$s\in S$ and $\1^*=\1\1=\1$ together with the structure on $S$ turns 
$S_\1$ into an involutive semigroup. 
\par\nin (b) Let ${\cal H}$ be a pre-Hilbert space. A {\it hermitian  
representation} $(\pi, {\cal H})$ of $S$ is a semigroup homomorphism 
$\pi\: S\to \End({\cal H})$ such that $\la \pi(s^*).v, w\ra= 
\la v, \pi(s).w\ra$ holds for all $s\in S$ and $v,w\in {\cal H}$. 
If in addition ${\cal H}$ is a Hilbert space and $\pi(S)\subeq 
B({\cal H})$, where $B({\cal H})$ denotes the bounded operators 
on ${\cal H}$, then we call $(\pi, {\cal H})$ a {\it representation}
of $S$.  
\par\nin (c) If $S$ is an involutive semigroup, then an {\it absolute 
value} on $S$ is a mapping $\alpha\: S\to [0,\infty[$ which satisfies 
$$ \alpha(s)=\alpha(s^*)\qquad\hbox{and}\qquad 
\alpha(st)\leq \alpha(s)\alpha(t)$$  
for $s,t\in S$. The collection of all absolute values on $S$ is
denoted by ${\cal A}(S)$.
\par\nin (d) A representation $(\pi, {\cal H})$ of $S$ is called
$\alpha$-{\it bounded} for some $\alpha\in {\cal A}(S)$ if $\|\pi(s)\|\leq 
\alpha(s)$ holds for all $s\in S$. \qed

\Definition I.2.  (Positive definite functions) Let $S$ be an
involutive semigroup. 
\par\nin (a) A function $\phi\: S\to \C$ is called {\it positive 
definite} if for all finite sequences $s_1, \ldots, s_n$ in $S$ and 
$c_1, \ldots, c_n$ in $\C$ one has 
$$ \sum_{j,k=1}^n c_j \oline{c_k}\phi(s_k s_j^*)\geq 0.$$ 
\par\nin (b) Let $\phi$ be a positive definite function on $S$. 
For each $t\in S$  we define a function $\phi_t(s)\:=\phi(st)$ on $S$ and set 
${\cal H}_\phi^0=\span\{ \phi_t\: t\in S\}$. Then the prescription 
$$\la \sum_{j=1} c_j \phi_{t_j}, \sum_{k=1} l_k \phi_{s_k}\ra
=\sum_{j,k} c_j \oline{l_k} \phi(s_k^*t_j)$$
defines a pre-Hilbert space structure on ${\cal H}_\phi^0$. We denote
by ${\cal H}_\phi$ the closure of ${\cal H}_\phi^0$ and note that 
${\cal H}_\phi$ also consists of functions on $S$ (cf.\ [Ne99, Ch.\ III]).
\par\nin (c) If $\phi$ is positive definite, then the prescription 
$$\pi_\phi^0\: S\to \End({\cal H}_\phi^0), \ \ (\pi_\phi^0(s).f)(t)=f(ts)$$
defines a hermitian representation $(\pi_\phi^0,
{\cal H}_\phi^0)$ of $S$ on the pre-Hilbert space ${\cal H}_\phi^0$. In the 
case where all operators  $\pi_\phi^0(s)$, $s\in S$, are bounded, 
the hermitian representation $(\pi_\phi^0, {\cal H}_\phi^0)$ 
extends to a representation $(\pi_\phi, {\cal H}_\phi)$ of $S$ which is
also given by right translations (cf.\ [Ne99, Th.\ III.1.3]). 
\par\nin (d) If $\alpha\in {\cal A}(S)$, then a positive 
definite function $\phi$ on $S$ is called $\alpha$-{\it bounded} if
$$(\exists C>0)(\forall s,t\in S)\qquad  |\phi(t^*st)|\leq 
C \alpha(s)\phi(t^*t).$$
Note that $C$ may be replaced by $C=1$ and that the
$\alpha$-boundedness
of $\phi$ is equivalent to the $\alpha$-boundedness of $(\pi_\phi,
{\cal H}_\phi)$ (cf.\ [Ne99, Th.\  III.1.19]). 
\par We denote by ${\cal P}(S,\alpha)$ the convex cone of all 
$\alpha$-bounded positive definite functions on $S$.
\qed

\Definition I.3. (Enveloping $C^*$-algebras) Let $S$ be an involutive 
semigroup and $\alpha$ be an absolute value on it. Then for a subset 
${\cal E}\subeq {\cal P}(S,\alpha)$ we set 
$${\cal H}_{\cal E}\:= \hat \bigoplus_{\phi\in {\cal E}}{\cal H}_\phi.$$ 
Then $\pi_{\cal E}(s)\:=\oplus_{\phi\in {\cal E}}\pi_\phi$ defines
an $\alpha$-bounded representation of $S$ on ${\cal H}_{\cal E}$. 
We denote by $C^*(S,\alpha, {\cal E})$ the closure of $\span (\pi_{\cal
E}(S))$ in $B({\cal H}_{\cal E})$ and note that this is a
$C^*$-subalgebra of $B({\cal H}_{\cal E})$. Note that the 
mapping 
$$j\: S\to C^*(S,\alpha, {\cal E}), \ \ s\mapsto \pi_{\cal E}(s)$$ 
has total image by definition. If ${\cal E}={\cal P}(S,\alpha)$, then
we  set $C^*(S,\alpha)\:=C^*(S,\alpha, {\cal E})$. \qed

\Lemma I.4. Every representation $(\tilde \pi, {\cal
H})$ of $C^*(S,\alpha)$ gives via 
$$\tilde \pi \mapsto \pi\:=\tilde\pi\circ j$$
rise to a $\alpha$-bounded representation $(\pi, {\cal H})$
of $S$. Moreover, this correspondence is bijective. 

\Proof. [Ne99, Th.\ III.2.9].\qed

Representations of involutive semigroups can be described in terms of
positive definite functions, those of $C^*$-algebras by positive
functionals. We describe now a correspondance between these two
pictures.

For every $C^*$-algebra ${\cal A}$ we denote by ${\cal A}_\1$ its
unification (adjoining of a unit element if there was none). Note that every positive 
functional $f$ on ${\cal A}$ extends uniquely to a positive functional 
$\tilde f$ on ${\cal A}_\1$ with $\tilde f(\1) = \|f\|$ (cf.\ [Dix82,
Sect.\ 2.4.3]).

\par We set 
$$E(C^*(S,\alpha)):=\{f\in C^*(S,\alpha)': f\ \hbox{positive
functional};
\  \|f\|=1\},$$
and 
$${\cal P}(S,\alpha)_\1\:=\{\phi\res_S\: \phi\in {\cal
P}(S_\1,\alpha), \ \phi(\1)=1\}.$$
Note that for all $\phi\in {\cal P}(S,\alpha)_\1$ one has
$\phi\in {\cal H}_\phi$ (cf.\ [Ne99, Prop.\ III.1.23]).

\Theorem I.5. If we equip ${\cal P}(S,\alpha)_\1$ with the
topology of pointwise convergence on $S$ and $E(C^*(S,\alpha))$ with
the weak-$*$-topology with respect to $C^*(S,\alpha)$, then the mapping  
$$\Psi\:{\cal P}(S,\alpha)_\1\to E(C^*(S,\alpha)),\ \ \phi\mapsto
f_\phi;\  f_\phi(x)\:=\la\tilde{\pi_\phi}(x).\phi, \phi\ra$$
is a homeomorphism.

\Proof. First we show that $\Psi$ is defined. Since $f_\phi$ is
obviously positive, we have to show that $\|f_\phi\|=1$ or equivalently
$\tilde f_\phi(\1)=1$. But since $(\tilde \pi_\phi, {\cal H}_\phi)$ is 
non-degenerate (cf.\ Lemma I.4),  it extends naturally to a representation of 
$C^*(S,\alpha)_\1$ by setting $\tilde {\pi_\phi}(\1)=\id$. Thus  

$$\tilde{f_\phi}(\1)=\la \tilde{\pi_\phi}(\1).\phi,\phi\ra=\la
\phi,\phi\ra=\phi(\1)=1,$$
and so $\Psi$ is defined. 
\par We describe now the construction of $\Psi^{-1}$. Let $f\in
E(C^*(S,\alpha))$. We consider $C^*(S,\alpha)$
as an involutive semigroup and $f$ as a positive definite function on
it. In this sense let $(\pi_f, {\cal H}_f)$ be the representation of 
Definition I.2(c). Since $f$ extends to a positive definite function
$\tilde f$ on $C^*(S,\alpha)_\1$, it follows from [Ne99, Prop.\ III.1.23] that $f\in 
{\cal H}_f$ and that $\pi_f$ extends to a representation of
$C^*(S,\alpha)_\1$ denoted by the same symbol. 
Note that $f(x)=\la\pi_f(x).f,f\ra$ for all $x\in C^*(S,\alpha)_\1$. Now we define 
$\phi_f\in {\cal P}(S,\alpha)$ by 
$\phi_f(s):=\la\pi_f(j(s)).f,f\ra$ for all $s\in S$.
{}From the fact that $f$ extends to a positive definite function on
$C^*(S,\alpha)_\1$, it follows that $\phi_f$ extends to a positive
definite function on $S_\1$ by setting $\phi_f(\1)=f(\1)=1$.

\par For the bijectivity of $\Psi$ we now have to show that 
$$(\forall \phi\in {\cal P}(S,\alpha)_\1)\quad
\phi=\phi_{f_\phi}\leqno(1.1)$$
$$(\forall f\in E(C^*(S,\alpha)) \quad  f=f_{\phi_f}.\leqno(1.2)$$
By the definition of $f_\phi$ we have 
$$f_\phi(j(s))=\la\tilde{\pi_\phi}(j(s)).\phi,\phi\ra=\phi(s)\leqno(1.3)$$
for all $s\in S$. Thus by the definition of $\phi_f$ we get that 
$$\phi_{f_\phi}(s)=f_\phi(j(s))=\phi(s)$$
for $s\in S$, proving  (1.1). 
By the same reasoning  we obtain that 
$$f_{\phi_f}(j(s))=\phi_f(s)=f(j(s))$$
so that $f$ and $f_{\phi_f}$ coincide on the total set $j(S)$. Since
both are continuous, they coincide, proving (1.2).

\par\nin Continuity of $\Psi$: If $\phi_n\to \phi$ pointwise on $S$,
then (1.3) implies that $f_{\phi_n}\to f_\phi$ on the dense set $\span\{j(S)\}$. Since 
$E(C^*(S,\alpha))$ is a bounded subset of continuous functions on
$C^*(S,\alpha)$ it thus follows that $f_{\phi_n}\to f_{\phi}$ in the weak-$*$-topology 
of $C^*(S,\alpha)$. 

\par\nin Continuity of $\Psi^{-\1}$: If $f_{\phi_n}\to f_\phi$
pointwise on $C^*(S,\alpha)$, we conclude from (1.3)  that 
$\phi_n\to \phi$ pointwise on $S$. \qed

\sectionheadline{II. The dual of complex Ol'shanski\u\i{} semigroups}

In this section we focus our interest on a concrete class of
involutive semigroups, namely complex Ol'shanski\u\i{}
semigroups, which may be thought as generalizations of  complex Lie
subsemigroups of complex Lie groups.

\subheadline {Complex Ol'shanski\u\i{} semigroups}

\Definition II.1. Let $\g$ be a finite dimensional real Lie algebra. 

\par\nin (a) An element $X\in \g$ is called {\it elliptic} if 
$\ad X$ is semsimple with purely imaginary spectrum. Accordingly we 
call a subset $W\subeq \g$ 
{\it elliptic} if all its elements are elliptic. 

\par\nin (b) A subalgebra $\a\subeq \g$ is said to be {\it compactly embedded}, 
if $\la e^{\ad\a}\ra$ is relatively compact in $\Aut(\g)$. Note that 
a subalgebra is compactly embedded if and only if it is elliptic. \qed

\Definition II.2. Let $\eset\neq W\subeq \g$ be an open convex $\Inn(\g)$-invariant elliptic cone and
$\oline W$ its closure.
Let $\tilde G$, resp.\ $\tilde G_\C$, be the simply connected 
Lie groups associated to $\g$, resp.\ $\g_\C$, and set 
$G_1\:=\la \exp\g\ra\subeq 
\tilde G_\C$. Then Lawson's Theorem (cf.\ [HiNe93, Th.\ 7.34, 35]) says that the 
subset $\Gamma_{G_1}(\oline W)\:=G_1\exp(i\oline W)$ is a closed subsemigroup of $G_\C$ 
and the polar map 
$$G_1\times \oline W\to \Gamma_{G_1}(\oline W), \ \ (g, X)\mapsto g\exp(iX)$$
is a homeomorphism. 
\par Now the universal covering semigroup $\Gamma_{\tilde G}(\oline W)\:=
\tilde\Gamma_{G_1}(\oline W)$ has a similar structure. We can lift the exponential function 
$\exp\: \g + i \oline W\to \Gamma_{G_1}(\oline W)$ to an exponential mapping 
$\Exp\: \g + i\oline W\to \Gamma_{\tilde G}(\oline W)$ with $\Exp(0)=\1$ and thus obtain 
a polar map $\tilde G\times \oline W\to \Gamma_{\tilde G} (\oline W),\ (g, X)\mapsto g\Exp(iX)$
which is a homeomorphism. 
\par If $G$ is  a connected Lie group associated to $\g$, then $\pi_1(G)$ is a 
discrete central subgroup of $\Gamma_{\tilde G}(\oline W)$ and we obtain a covering 
homomorphism $\Gamma_{\tilde G}(\oline W)\to \Gamma_G(\oline
W)\:=\Gamma_{\tilde G}(
\oline W)/ \pi_1(G)$  
(cf.\ [HiNe93, Ch.\ 3]). It is easy to see that there is also a polar map 
$G\times \oline W \to \Gamma_G(\oline W), (g, X)\mapsto g\Exp(iX)$ which is a homeomorphism. 
The semigroups of the type $\Gamma_G(\oline W)$ are called {\it complex Ol'shanski\u\i{}
semigroups}.
\par The subset $\Gamma_G(W)\subeq \Gamma_G(\oline W)$ is an open semigroup carrying 
a complex manifold structure such that the multiplication is holomorphic. 
Moreover there is an involution on $\Gamma_G(\oline W)$ given by 

$$^*\: \Gamma_G(\oline W)\to \Gamma_G(\oline W), s=g\Exp(iX)\mapsto s^*=\Exp(iX)g^{-1}$$
being antiholomorphic on $\Gamma_G(W)$ (cf.\ [HiNe93, Th.\ 9.15] for a
proof of all that).
Thus both $\Gamma_G( W)$ and $\Gamma_G(\oline W)$  are involutive 
semigroups.\qed

{}From now on we denote by $S$ an open complex Ol'shanski\u\i{}
semigroup $\Gamma_G(W)$ and by $\oline S$ its ``closure'' $\Gamma_G(\oline
W)$.

\Definition II.3. A {\it holomorphic representation} $(\pi,{\cal H})$ 
of a complex Ol'shanski\u\i{} semigroup $S$ is a non-degenerate 
representation of the involutive semigroup $S$ (cf.\ Definition I.1(b))  
for which the map $\pi\: S\to B({\cal H})$ is 
holomorphic. \qed 

Since we now work in the category of complex Ol'shanski\u\i{} 
semigroups and holomorphic representations we have to adapt 
the objects ${\cal A}(S)$  and ${\cal P}(S,\alpha)$ slightly to the new
setting.

\Definition II.4.  Let $V$ be a finite dimensional real vector space 
and $V^*$ its dual space. 
\par\nin (a) For a subset $E\subeq V$ we define its {\it dual set} by 
$E^\star\:=\{ \alpha\in V^*\:(\forall X\in E) \alpha(X)\geq 0\}$. Note
that $E^\star$ is a closed convex cone in $V^*$.    
\par\nin (b) If $C\subeq V$ is a convex set, then we define 
$$\lim C\:=\{ v\in V\: v+C\subeq C\}\qquad \hbox{and}\qquad B(C)\:=\{
\alpha\in V^*\: \inf \alpha (C)>-\infty\}$$ 
We call $\lim C$ the {\it limit cone} of $C$. 
Note that both $\lim C$ and $B(C)$ are convex sets and that $\lim C$ 
is closed if $C$ is open or closed. 
\par\nin (c) If $C\subeq V$ is a convex set, then we say $C$ is
{\it pointed} if $C$ contains no affine lines, i.e., $\oline{\lim C}\cap -\oline
{\lim C}=\{0\}$. \qed

\Lemma II.5. Let $(\pi, {\cal H})$ be a holomorphic representation of
$S$. Then $\|\alpha(s)\|\:=\|\pi(s)\|$ is a continuous $G\times
G$-biinvariant absolute value on $S$ and there exists a unique closed  
$\Ad(G)^*$-invariant convex subset $C\subeq \g^*$ with 
$-W^\star =\lim C $  such that 
$$\alpha(g\Exp(iX))=\: \alpha_C(g\Exp(iX))\:= e^{\sup\la X,C\ra}$$
for $g\Exp(iX)\in S$. 
\Proof. [Ne99, Th.\ XI.3.5].\qed

Note that the condition $-W^\star=\lim C$ is equivalent to 
$-\oline W= \oline {B(C)}$ (cf.\ [Ne99, Lemma V.1.18]).
We define now 
$${\cal A}_h(S)\:=\{ \alpha_C\: C\subeq \g^*,\ \hbox{closed and convex},
\ -W^\star =\lim C, \Ad(G)^*.C=C\},$$ 
and for all $\alpha\in {\cal A}_h(S)$ we set 
$${\cal P}_h(S,\alpha)\:=\{\phi\in \Hol(S)\: \phi\ \hbox{positive
definite and $\alpha$-bounded}\}.$$
Finally we set $C_h^*(S,\alpha)\:=C^*(S,\alpha, {\cal P}_h(S,\alpha))$
(cf.\ Definition I.3). 

Since each $\alpha$ is locally bounded, the canonical mapping 
$j\: S\to C_h^*(S,\alpha)$ becomes holomorphic. The analog of Lemma
I.4 now reads as follows.

\Lemma II.6. The prescription $\tilde \pi\mapsto \pi\:=\tilde\pi\circ
j$ defines a bijection between non-degenerate representations 
$(\tilde \pi, {\cal H})$ of $C_h^*(S,\alpha)$ and $\alpha$-bounded 
holomorphic representations $(\pi, {\cal H})$ of $S$.\qed

\Definition II.7.  A $C^*$-algebra ${\cal A}$ is said to be $CCR$ ({\it completely
continuous representations}) or {\it liminal} if for all irreducible representations
$(\pi,{\cal H})$ of ${\cal A}$ the image $\pi({\cal A})$ is contained in the
space of compact operators ${\cal K}({\cal H})$ of ${\cal H}$.\qed

\Theorem II.8. {\rm (K.-H. Neeb)} If $S$ is a
complex Ol'shanski\u\i{} semigroup and $\alpha\in {\cal A}_h(S)$, then 
the $C^*$-algebra $C_h^*(S,\alpha)$ is CCR.

\Proof. [Ne99, Th.\ XI.6.6].\qed

\subheadline{The hull-kernel topology on the dual}

In this subsection we introduce for each $\alpha\in {\cal A}_h(S)$ the
hull-kernel topology on the set $\hat S_\alpha$ of equivalence classes
of irreducible $\alpha$-bounded holomorphic representations of $S$. 
Then we characterize the hull-kernel topology on 
$\hat S_\alpha$ in terms of the compact convergence of the corresponding  matrix coefficients 
of irreducible holomorphic representations.

\Lemma II.9. On ${\cal P}_h (S,\alpha)$ the topology of pointwise
convergence coincides with the topology of compact convergence.

\Proof. Let $\phi_n\to\phi$ pointwise on $S$ and $K\subeq S$ a compact
subset. Then we find an element $t\in S$  and a compact subset $\tilde
K\subeq S$ such that $K\subset t^*\tilde K t$ holds (cf.\  [HiNe93, 3.19]).
Since for all $\psi\in{\cal P}_h(S,\alpha)$ and $s,t\in S$ the
inequality 

$$|\psi(t^*st)|\leq \alpha(s)\psi(t^*t)$$
holds (cf.\ Definition I.2(d)), the locally boundedness of $\alpha$
implies that all functions in ${\cal P}_h(S,\alpha)$ are uniformly
bounded on $K$. Thus by Montel's Theorem ${\cal P}_h(S,\alpha)\res_K$ is a
compact subset of $C(K)$. In particular, $\phi_n\to \phi$ pointwise on
$S$ implies $\phi_n\to \phi$ uniformly on
$K$, as was to be shown. \qed

\Theorem II.10. If we equip ${\cal P}_h(S,\alpha)_\1$ with the
topology of compact convergence on $S$ and $E(C_h^*(S,\alpha))$ with the weak-$*$-topology of 
$C_h^*(S,\alpha)$, then the mapping 

$$\Psi\: {\cal P}_h(S,\alpha)_\1\to E(C_h^*(S,\alpha)),\ \ \phi\mapsto f_\phi$$
is a homeomorphism.

\Proof. In view of Lemma II.6 and Lemma II.9, this follows now from 
Theorem I.5. \qed

\Definition II.11. (Hull-kernel topology) Let  ${\cal A}$ be a
$C^*$-algebra. 
\par\nin (a) A two-sided ideal $J$ of ${\cal A}$ is called {\it
primitive} if there exists an irreducible representation $(\pi,{\cal
H})$ of ${\cal A}$ such that $J=\ker \pi$. The space of primitive
ideals is denoted by $\Prim ({\cal A})$. 
\par\nin (b) For a subset $M\subeq \Prim({\cal A})$ we set 
$$I(M):=\bigcap_{J\in M}J\qquad\hbox{and}\qquad\oline M:=\{J\in\Prim
({\cal A}):J\supset I(M)\}.$$ 
Then the prescription 
$${\cal T}_J:=\{\Prim ({\cal A})\bs \oline M:M\subset\Prim({\cal A})\}$$   
defines a topology on  $\Prim({\cal A})$ (cf.\ [Dix82, Sect.\ 3.1]), the 
so-called {\it Jacobson topology}. 
\par\nin (c) Let $\hat {\cal A}$ denote the set of equivalence classes
of irreducible representations of ${\cal A}$. For each irreducible 
representation $(\pi,{\cal H}_\pi)$ of ${\cal A}$ we denote by $[\pi]$ its 
equivalence class. Consider the map 
$$\hat{\cal A}\to\Prim({\cal A}),\ [\pi]\mapsto\ker\pi.$$
If we equip $\Prim({\cal A})$ with the Jacobson topology, then the
initial topology on $\hat{\cal A}$ under this map is called the {\it
hull-kernel} or {\it Fell topology} on $\hat{\cal A}$. \qed

\Definition II.12. Let $S$ be a complex Ol'shanski\u\i{} semigroup
and $\alpha\in {\cal A}_h(S)$ an absolute value on it. 
We denote by $\hat S$ the set of equivalence classes of 
irreducible holomorphic representations of $S$ and by $\hat S_\alpha$
the subset of all $\alpha$-bounded ones. Note that $
C_h^*(S,\alpha)\,\hat{}$ can be identified with $\hat S_\alpha$ (cf.\
Lemma II.6). The topology on $
\hat S_\alpha$ induced from the hull-kernel topology on
$C_h^*(S,\alpha)\, \hat{}$
is denoted by ${\cal T}_{hk}^\alpha$. Note that $\hat S_\alpha$
is countable at infinity, since $C_h^*(S,\alpha)$ is separable. \qed

\Remark II.13.  (a) If $(\pi, {\cal H}_\pi)$ is an $\alpha$-bounded  holomorphic 
representation of 
$S$ and $v\in {\cal H}_\pi$, then an element 
of ${\cal P}_h(S,\alpha)_\1$ is defined by 
$\pi_{v,v}(s)\:=\la \pi(s).v, v\ra$, $s\in S$. 
Conversely, the Gelfand-Naimark-Segal correspondance between 
positive definite functions and matrix coefficients together 
with [Ne99, Prop.\ III.1.23] imply that 
$${\cal P}_h(S,\alpha)_\1=\{ \pi_{v,v}\: v\in
{\cal H}_\pi,\  (\pi, {\cal H}_\pi)\ \alpha-\hbox{bounded}  \}.$$
\par\nin (b) Every holomorphic representation 
$(\pi, {\cal H})$ of $S$
gives rise to a representation of the involutive semigroup
$S\cup G$, also denoted by $(\pi, {\cal H})$, 
such that $\pi\res_G$ is a uniquely
determined unitary representation of $G$. 
Moreover, $\pi\res_G$ 
is irreducible if and only if $\pi\res_S$ is (cf.\ [Ne99, Ch.\ XI] for 
all that). 
\par\nin (c) In view of (a) and (b), we conclude that every $\phi\in 
{\cal P}_h(S,\alpha)_\1$ extends naturally to a positive definite 
function on $S\cup G$ which is continuous when restricted to $G$. 
In the sequel we identify the elements of ${\cal P}_h(S,\alpha)_\1$
as positive definite functions on $S\cup G$. Further 
for all  $\phi\in {\cal P}_h(S,\alpha)_\1$
all the {\it radial limits} 
$$\lim_{t\to 0\atop t>0}\phi(g\Exp(itX))=\phi(g)$$ 
exist, where $g\in G$ and $X\in W$ (this is immediate from (a) and 
[Ne99, Prop.\ IV.3.2]). \qed

\Theorem II.14. For a subset $M$  of $\hat S_\alpha$ let $\oline M$ be
the closure in the topology ${\cal T}_{hk}^\alpha$ (cf.\ {\rm
Definition II.12}).
For $[\pi]\in\hat S_\alpha$ and $v\in {\cal H}_\pi$ let 
$\pi_{v,v}(s):=\la \pi(s).v,v\ra$,
$$E([\pi]):=\{\pi_{v,v}: v\in{\cal H}_\pi,\ \|v\|=1\}\ \
\hbox{and}\ \ E(M):=\bigcup_{[\pi]\in M}E([\pi]).$$
Let $\oline{E(M)}$ be the closure of $E(M)$  in ${\cal P}_h(S,\alpha)_\1$,
equipped with the topology of compact convergence on $S$.
Then the following assertions are equivalent:

\item{(1)}$[\pi]\in\oline M$,
\item{(2)}$(\exists v\in{\cal H}_\pi,\ \|v\|=1)\quad \pi_{v,v}\in\oline{E(M)}$,
\item{(3)}$(\forall v\in{\cal H}_\pi,\ \|v\|=1)\quad \pi_{v,v}\in\oline{E(M)}$,
\item{(4)}$E([\pi])\cap\oline{\span E(M)}\neq\{0\}$.

\Proof. In view of Theorem II.10 and Remark II.13(a),  
the assertion follows now from 
[Dix82, Th.\ 3.4.10].\qed

\subheadline{Comparison with the topology on $\hat G$}

Recall from Remark II.13(b) that every irreducible holomorphic 
representation of $S$
defines in a unique manner an irreducible unitary representation of
$G$. Thus we may consider $\hat S$ as a certain subset of $\hat G$ and
the topology on $\hat G$ gives rise to a topology on $\hat S_\alpha$
which we denote by ${\cal T}_G^\alpha$. 
\par In this subsection we compare the hull-kernel topology 
${\cal T}_{hk}^\alpha$ on $\hat S_\alpha$ with the topology 
${\cal T}_G^\alpha$ induced from $\hat G$. With real
analysis methods, i.e., boundary values of Poisson transforms, we show
that ${\cal T}_{hk}^\alpha\subeq {\cal T}_G^\alpha$.

\Lemma II.15. Let $\alpha\in {\cal A}_h(S)$. Then for each compact subset
$Q\subeq W$ there exists a norm $\|\cdot\|$ on $\g$ and a constant 
$c>0$ sucht that 

$$(\forall X\in \R^+ Q) \qquad \alpha(\Exp(iX))\leq ce^{\|X\|}.$$

\Proof. If $0\in W$, then $W=\g$ and $\alpha =1$ by Lemma II.5. In
this case the assertion of the lemma is clear. Thus we may assume that $0\not\in
W$.
\par Note that $\alpha =\alpha_C$ for some closed convex subset
$C\subeq \g^*$ with $-\oline W=\oline {B(C)}$ by the definition of 
${\cal A}_h(S)$. Thus the compactness of $Q$ in
$W$ shows that 
$$\sup_{X\in [0,1]Q} \alpha (\Exp(iX))=\:c<\infty.$$  
Let $\|\cdot\|$ be a norm on $\g$ such that $\|X\|>1$ for all $X\in 
]1, \infty[Q$ (this is possible since $0\not\in W$). 
Note that the mapping $W\to \R^+, \ 
X\mapsto \log\alpha(\Exp(iX))$ is positively homogeneous by the
construction of ${\cal A}_h(S)$. Thus is $Y=kX\in \R^+ Q$ with $k\in
\R^+$ and $X\in Q$, then  
$$\alpha(\Exp(iY))=\alpha(\Exp(iX))^k\leq c^k\leq c^{\|Y\|}.$$
Now an appropriate rescaling of $\|\cdot\|$ yields the assertion.\qed

\Proposition II.16. The mapping  $(\hat S_\alpha,{\cal T}_G^\alpha)\to (\hat S_\alpha,
{\cal T}_{hk}^\alpha)$ is continuous, i.e., we have ${\cal
T}_{hk}^\alpha\subset
{\cal T}_G^\alpha$.

\Proof. Recall from Remark II.13(c) that the elements of  
${\cal P}_h(S,\alpha)_\1$ are
identified with the positive definite functions on $S\cup G$
which restrictions $\phi\res_S$ are holomorphic and 
$\alpha$-bounded and have $\phi\res_G$ as continuous radial 
limit. 
\par In view of  [Dix82, Prop.\ 18.1.5] and Theorem II.14 it suffices  to
show that  a sequence
$(\phi_n)_{n\in\N}$ in ${\cal P}_h(S,\alpha)_\1$ 
which satisfies 
$$\phi_n\to\phi\ \ \hbox{ compactly on $G$},$$ 
also satisfies  
$$\phi_n\to\phi\ \ \hbox{compactly on  $S$}.$$
\par Let $K\subeq S$ be a compact subset and  $s=g\Exp(iX)\in K$.
Let $X_1,\ldots,X_n$ be a basis of $\g$ which is contained in  $W$ and
set $Q:=\conv(\{X_1,\ldots,X_n\})$. We choose the basis in such a way 
that  $X\in Q$. Since 
$d\Exp(iX)$ is everywhere invertible, we find an open subset 
$U\subset \g+iQ$ with  $iX\in U$ such that  $\Exp\res_U$ is a 
diffeomorphism.
W.l.o.g. we may assume $U\subeq Q$ and also $K=g\Exp(U)$.
Let $\R_+\:=]0,\infty[$. Then the mapping 
$$h\:\Gamma_{\R^n}(\R_+^n)\to S,\ \ (z_1,\ldots,z_n)\mapsto 
g\Exp(z_1X_1+\ldots +z_nX_n),$$
induces a map 
$$ h_*:\Hol(S)\to \Hol(\Gamma_{\R^n}(\R_+^n)),\ \
\psi\mapsto\tilde\psi:=\psi\circ h. $$
\par Thus we only have to show that  
the compact convergence of $\tilde{\phi_n}\res_{\R^n}$ implies the
compact convergence of $\tilde{\phi_n}$ on the tube domain 
$\Gamma_{\R^n}(\R_+^n)$.  Recall that 
each $\psi\in {\cal P}_h(S,\alpha)_\1$ 
satisfies the estimate 
$$(\forall s\in S)\quad |\psi(s)|\leq \psi(\1) \alpha(s)=\alpha(s)$$
(cf.\ Definition I.2(d)).  Let $c>0$ and define $f\in
\Hol(\Gamma_{\R^n}(\R_+^n))$
$$f:\Gamma_{\R^n}(\R_+^n)\to\C,\ \ f(z):=ce^{-ic(z_1+\ldots+z_n)}.$$ 
In view of Lemma II.15, we then have $|\tilde \psi(z)|\leq |f(z)|$ for
all $z\in \Gamma_{\R^n}(\R_+^n)$ and $\psi\in {\cal
P}_h(S,\alpha)_\1$ provided $c>0$ is chosen sufficiently large. 
Since $f$ has no zeros, we may replace $\tilde\phi_n$ by 
$\phi_n':={1\over f}\tilde\phi_n$. 
Note that the functions $\phi_n'$, $\phi'$  are
elements of $\Hol(\Gamma_{\R^n}(\R_+^n))$ uniformly bounded by
$1$. Further we replace $\phi_n'$, $\phi$  by $\phi_n'':= g\phi_n'$
and $\phi'':= g\phi'$, where 
$g\:\Gamma_{\R^n}(\R_+^n)\to\C,\ g(z):={1\over z_1+i}\cdot\ldots\cdot
{1\over z_n+i}$. Note that the functions 
$\phi_n''$, $\phi''$ are uniformly bounded 
elements of $\Hol(\Gamma_{\R^n}(\R_+^n))$ which are uniformly 
vanishing at infinity. 
\par Let 
$$p\:\oline{\Gamma_{\R^n}(\R_+^n)}\to \R,\ \ p(x_1+iy_1,\ldots,x_n+iy_n):=\Big({1\over \pi}\Big)^n
\prod_{j=1}^n{y_j\over x_j^2+y_j^2}$$
be the Poisson kernel of the upper half plane and 
$$P\:L^\infty(\R^n)\to{\rm Harm}(\Gamma_{\R^n}(\R_+^n)),\ \ P(f)(z):=\int_{\R^n}
f(x)p(z-x)\ dx$$ 
the corresponding Poisson transform.
According to [SW75, \S2, 2.1,5], the functions $\phi_n''$, $\phi''$ are the 
Poisson transforms of their boundary values, 
i.e.,  $\phi_n''=P(\phi_n''\res_{\R^n})$, 
$\phi''=P(\phi''\res_{\R^n})$ and we have 
$$\|\phi_n''\|_{\Gamma_{\R^n}(\R^n_+),\infty}=\|\phi_n''\|_{\R^n,\infty}\leqno(2.1)
$$
and analogously for $\phi''$. 
Since the functions $\phi_n''\res_{\R^n}$, $\phi''\res_{\R^n}$ vanish uniformly at infinity, 
equation (2.1) implies that the compact convergence of the
$\phi_n''\res_{\R^n}\to \phi''\res_{\R^n}$ implies the 
compact convergence on $\Gamma_{\R^n}(\R_+^n)$. This concludes the
proof of the proposition. \qed

\subheadline{Highest weight representations}

In this section we take a closer look at the irreducible holomorphic
representations  of $S$. It turns out that  they are obtained by
analytic continuation of unitary highest weight representations of
$G$. 

\par Note that if a real Lie algebra admits a non-empty open elliptic convex
cone, then there exists a compactly embedded Cartan subalgebra
$\t\subeq \g$ (cf.\ [Ne99, Th.\ VII.1.8]). To step further we first need some
terminology concerning Lie algebras with compactly embedded Cartan
subalgebras.

\Definition II.17. Let $\g$ be a finite dimensional Lie algebra over
$\R$ with compactly embedded Cartan subalgebra $\t$.

\par\nin (a) Associated to the Cartan subalgebra $\t_\C$ in the complexification 
$\g_\C$ there is a root decomposition as follows. For a linear functional 
$\alpha\in \t_\C^*$ we set 
$$\g_\C^\alpha\:=\{X\in \g_\C\: (\forall Y\in \t_\C)\ [Y,X]=\alpha(Y)X\}$$  
and write $\Delta\:=\{\alpha\in\t_\C^*\bs\{0\}\:\g_\C^\alpha\neq \{0\}\}$ for the 
set of roots. Then $\g_\C=\t_\C\oplus\bigoplus_{\alpha\in \Delta}\g_\C^\alpha$, 
$\alpha(\t)\subeq i\R$ for all $\alpha\in \Delta$ and $\oline{\g_\C^\alpha}=
\g_\C^{-\alpha}$, where $X\to\oline X$ denotes complex conjugation on $\g_\C$
with respect to $\g$. 

\par\nin (b)  Let $\k$ be a maximal compactly embedded
subalgebra of $\g$ containing $\t$. Then a root $\alpha$ is said to be
{\it compact} 
if $\g_\C^\alpha\subeq \k_\C$ and {\it non-compact} otherwise. We write $\Delta_k$for the set of 
compact roots and $\Delta_n$ for the non-compact ones.

\par\nin (c) A positive system $\Delta^+$ of roots is a subset of $\Delta$ for which there 
exists  a regular element $X_0\in i\t^*$ with $\Delta^+\:=\{ \alpha\in 
\Delta\: \alpha(X_0)>0\}$. We call a positive system $\k${\it-adapted} if the set 
$\Delta_n^+\:=\Delta_n\cap\Delta^+$ is invariant under the {\it Weyl group}
${\cal W}_\k\:=N_{\Inn(\k)}(\t)/ Z_{\Inn(\k)}(\t)$ acting on $\t$. We recall from 
[Ne99, Prop.\ VII.2.14] that there exists a $\k$-adapted positive system if and only if 
$\z_\g(\z(\k))=\k$. In this case we say $\g$ is {\it quasihermitian}. In this case it is 
easy to see 
that $\s$ is quasihermitian, too, and so all simple ideals of $\s$ are either 
compact or hermitian.  
 
\par\nin (d) We associate to the positive system $\Delta^+$ the convex cones 

$$C_{\rm min}\:=\oline{\cone\{ i[\oline {X_\alpha},X_\alpha]\: X_\alpha\in \g_\C^\alpha, 
\alpha\in \Delta_n^+\}},$$
and $C_{\rm max}\:=(i\Delta_n^+)^\star=\{X\in\t\: (\forall\alpha\in \Delta_n^+)
\ i\alpha(X)\geq 0\}$. Note that both $C_{\rm min}$ and $C_{\rm max}$ are closed 
convex cones in $\t$. 

\par\nin (e) Write $p_\t\:\g\to \t$ for the orthogonal projection along $[\t,\g]$ 
and set ${\cal O}_X\:=\Inn(\g).X$ for the adjoint orbit through $X\in \g$.
We define the {\it minimal} and {\it maximal cone} associated to $\Delta^+$ by 

$$W_{\rm min}\:=\{X\in \g\: p_\t({\cal O}_X)\subeq C_{\rm min}\}\quad\hbox{and}
\quad W_{\rm max}\:=\{ X\in \g\: p_\t({\cal O}_X)\subeq C_{\rm max}\}$$   
and note that both cones are convex, closed and $\Inn(\g)$-invariant. \qed

{}From now on we assume that $\g$ contains a compactly embedded Cartan subalgebra 
$\t\subeq \g$ and that there exists a non-empty  open elliptic  convex
cone $W\subeq \g$. Then in view of [Ne99, Th.\ VII.3.8],  there 
exits a $\k$-adapted positive system $\Delta^+$ such that  
$$W_{\rm min}\subeq \oline W\subeq W_{\rm max}$$
holds, $W_{\rm max}^0$ is elliptic, $W_{\rm min}\cap \t= C_{\rm min}$
and $W_{\rm max }\cap\t= C_{\rm max}$. Recall that every 
elliptic $\Ad(G)$-invariant cone $W\subeq \g$ can be 
reconstructed by its intersection with $\t$, i.e., $W=\Ad(G).(W\cap\t)$.

\Definition II.18.  Let $\Delta^+$ be a positive system. 

\par\nin (a) For a $\g_{\C}$-module $V$ and $\beta \in (\t_{\C})^*$
we write $V^\beta := \{ v \in V : (\forall X \in \t_{\C}) 
X.v = \beta(X)v \}$ for the 
{\it weight space of weight} $\beta$ and 
${\cal P}_V = \{ \beta \: V^\beta \not= \{0\} \}$ 
for the set of weights of $V$. 

\par\nin (b) Let $V$ be a $\g_{\C}$-module and $v \in V^\lambda$ a
$\t_\C$-weight vector. 
We say that $v$ is a {\it primitive element of V} (with respect to 
$\Delta^+$) 
if $\g_\C^\alpha.v = \{0\}$ holds for all $\alpha \in \Delta^+$. 

\par\nin (c) A $\g_{\C}$-module $V$ is called a {\it highest weight 
module} with highest weight $\lambda$ (with respect to $\Delta^+$) 
if it is generated by a primitive element of weight $\lambda$. 

\par\nin (d) Let $G$ be a connected Lie group with Lie algebra $\g$. 
We write $K$ for the analytic subgroup of $G$ corresponding to $\k$. 
Let $(\pi, {\cal H})$ be a unitary representation of $G$. A vector 
$v\in{\cal H}$ is called $K${\it -finite} if it is contained in a
finite dimensional $K$-invariant subspace. 
We write ${\cal H}^{K,\omega}$ for the space of analytic 
$K$-finite vectors. 

\par\nin (e) An irreducible unitary representation 
$(\pi_\lambda,{\cal H}_\lambda)$ of 
$G$ is called a {\it highest weight representation} with respect to $\Delta^+$ 
and highest weight $\lambda\in i\t^*$ if 
$L(\lambda)\:={\cal H}_\lambda^{K,\omega}$ is 
a highest weight module for $\g_\C$ with respect to $\Delta^+$ and 
highest weight $\lambda$. 
We write $HW(G, \Delta^+)\subset i\t^*$ for the set of highest weights 
corresponding to unitary highest weight representations of $G$ with
respect to $\Delta^+$.\qed

\Lemma II.19. Let $S=\Gamma_G(W)$ be a complex Ol'shanski\u\i{}
semigroup and $\Delta^+$ be a $\k$-adapted positive system with 
$C_{\rm min}\subeq \oline W\cap\t\subeq C_{\rm max}$.
\item{(i)} If $(\pi, {\cal H})$ is an irreducible holomorphic
representation of $S$, then $(\pi, {\cal H})$ extends to a 
representation of $S\cup G$, also denoted by 
$(\pi, {\cal H})$, such that $\pi\res_G$
is a uniquely determined unitary highest weight representation of $G$ 
with respect to $\Delta^+$. 
Conversely, if $(\pi_\lambda, {\cal H}_\lambda)$ is a unitary highest
weight representation of $G$ with respect to  $\Delta^+$, then 
$(\pi_\lambda, {\cal H}_\lambda)$
extends to a uniquely determined holomorphic representation of
$S$. In particular, the mapping 
$$HW(G, \Delta^+)\to \hat S, \ \ \lambda\mapsto [\pi_\lambda]$$
is bijective, hence gives a parametrization of $\hat S$ by highest
weights. 

\item{(ii)} If $\alpha=\alpha_C$ and $HW_\alpha$ denotes the subset of
$HW(G,\Delta^+)$ corresponding to $\alpha$-bounded representations, then 
$$HW_\alpha=\{ \lambda\in HW(G, \Delta^+)\: i\lambda\in
C\cap\t^*\}.$$

\Proof. (i) This follows from [Ne99, Th.\ XI.2.3].
\par\nin (ii) It follows from [Kr99a, Lemma IV.12] that 
$$HW_\alpha=\{ \lambda\in HW(G, \Delta^+)\: (\forall X\in W\cap\t)
 \ e^{\lambda(iX)}\leq \alpha(\Exp(iX))\}.$$
We define $f\: \t\to \R\cup\{+\infty\},\ X\mapsto \sup\la X, C\ra$ and
note that $f$ is a lower semicontinuous convex function.  Since $C$ was supposed  to
be closed, convex and $\Ad^*(G)$-invariant, it follows in particular
that $f(X)=\sup\la X, C\cap\t^*\ra$ for all $X\in \t$ (cf.\
 [Ne99,Prop.\ V.2.2]). 
\par Fix $\lambda\in HW(G, \Delta^+)$. Then we have 
$\lambda\in HW_\alpha$ if and only if $i\lambda(X)\leq f(X)$ for all 
$X\in W\cap\t$. Let $D_f\:=\{ X\in \t\: f(X)<\infty\}$. Then 
$\oline W= -\oline {B(C)}$ implies that 
$D_f=(-B(C))\cap\t\subeq \oline W\cap\t=\oline {W\cap\t}$. Thus it
follows from [Ne99, Prop.\ V.3.2(i)] that $\lambda\in HW_\alpha$ if
and only if $i\lambda\leq f$ on $\t$, i.e., $i\lambda\in C\cap\t^*$ as
was to be shown.\qed

\subheadline{ The topology on $\hat S_\alpha$}

In view of Lemma II.19, we can parametrize $\hat
S_\alpha$ by the subset $HW_\alpha\subeq i\t^*$.  The euclidean
topology on $HW_\alpha$ thus gives rise to a topology on $\hat
S_\alpha$ which we denote by ${\cal T}_e^\alpha$.

In the sequel we will show that for generic absolute values $\alpha$ 
all the topologies ${\cal T}_{hk}^\alpha$, 
${\cal T}_G^\alpha$ and ${\cal T}_e^\alpha$
coincide on $\hat S_\alpha$. 

\par Associated to the complex Lie subalgebras 
$\k_\C$, $\p_\C^+:=\bigoplus_{\alpha\in \Delta_n^+}
\g_\C^\alpha$ and $\p_\C^-:=\bigoplus_{\alpha\in -\Delta_n^+}\g_\C^\alpha$ of
$\g_\C$ we have analytic subgroups $K_\C$, $P^+$ and
$P^-$ of $G_\C$. Recall from [Ne99, Ch.\ XII] that the groups $P^\pm$ are
biholomorphic to $\p^\pm$ via the exponential mapping and that the
multiplication mapping 
$$P^-\times K_\C\times P^+\to G_\C,\ \ (p_-,k,p_+)\mapsto p_- kp_+$$
is biholomorphic onto its open image  $P^-K_\C P^+\subeq G_\C$. 
Set $S_1\:=\Gamma_{G_1}(W)\subeq G_\C$ and recall from [Ne99, Th.\ XII.4.6] that

$$G_1\subset \oline{S_1}\subset P^-K_\C P^+.$$
Further, if $P^-\tilde{K_\C} P^+\cong P^-\times \tilde{K_\C}\times
P^+$ denotes the simply connected covering of $P^- K_\C P^+$, then the
chain of inclusions from above lifts to
$$\tilde G\subset\oline{\tilde S}\subset P^-\tilde{K_\C} P^+$$
(cf.\ [Ne99, Cor.\ XII.4.7]). We denote by  
$$\kappa\: P^- \tilde {K_\C} P^+ \to\tilde{K_\C}, \ \ s\mapsto
 \kappa(s)$$
the holomorphic projection on the middle component. 

\par If $(\pi_\lambda,
{\cal H}_\lambda)$ is a unitary highest weight representation of $G$,
then we denote by $v_\lambda\in {\cal H}_\lambda$
a normalized highest weight vector. Further we set  
$F(\lambda)\:=U(\k_\C).v_\lambda$ for the minimal $\k$-type and write $(\pi_\lambda^K,
F(\lambda))$ for the finite dimensional holomorphic representation of 
$\tilde{K_\C}$ with highest weight $\lambda$. 
Finally we denote by $HW(\Delta_k^+)$ the set of linear functionals on
$\t_\C^*$ which are dominant integral with respect to $\Delta_k^+$
and note that $HW(G,\Delta^+)\subeq HW(\Delta_k^+)$.

\Proposition II.20. The mapping 
$$\Phi\: S\times HW(G,\Delta^+)\to\C,\ \ (s,\lambda)\mapsto\psi_\lambda(s):=\la 
\pi_\lambda(s).v_\lambda,v_\lambda\ra$$
lifts to a continuous mapping $\tilde \Phi\: P^-\tilde{K_\C}P^+\times
HW(\Delta_k^+)\to\C$ which is holomorphic in the first variable and
given explicitly by $\tilde \Phi(s,\lambda)=\la 
\pi_\lambda^K(\kappa(s)).v_\lambda,v_\lambda\ra$.

\Proof. W.l.o.g. we may assume that $S=\tilde S$. 
Fix $\lambda\in HW(G,\Delta^+)$.  Since $v_\lambda$ is an analytic
vector of $(\pi_\lambda, {\cal H}_\lambda)$, we find a zero-neighborhood  
$U\subset\g$ such that on $U_\C:=U+iU\subset\g_\C$
the series 
$$\sum_{j=0}^\infty{1\over j!}\la d\pi_\lambda(X)^j.v_\lambda,v_\lambda\ra
\leqno(2.2)$$
converges uniformly (cf.\ [Ne99, Lemma XI.2.1]). Note that the value of (2.2) coincides
with $\psi_\lambda(\exp(X))$ for all $X\in U$. 

\par\nin Let $D\subeq \g_\C$ be the connected component of $\{0\}$ in 
$\exp_{G_\C}^{-1}(P^-K_\C P^+)$ and $\Exp\:D \to P^-\tilde
K_\C P^+$ the lifting of $\exp_{G_\C}\res_D$ satisfying  $\Exp(0)=\1$.
If we choose $U$ sufficiently small, we may assume that 
$\Exp(U_\C)\subset P^-\tilde{ K_\C} P^+$ and that  $\Exp\res_{U_\C}$ is
a diffeomorphism onto its image. 
Let $V^-$, $V_0$, $V^+$ in $P^-$, $\tilde {K_\C}$, $P^+$ be connected
open $\1$-neighborhoods with $V^-V_0V^+\subset \Exp(U_\C)$.
\par\nin According to the the uniform convergence of (2.2) on $U_\C$, 
the prescription 
$$\psi'_\lambda\:\Exp(U_\C)\to\C, \ \ \psi'_\lambda(\Exp(X)):=
\sum_{j=0}^\infty{1\over j!}\la d\pi_\lambda(X)^j.v_\lambda,v_\lambda\ra$$
defines a holomorphic function on $\Exp(U_\C)$. 
We claim that 
$$\psi_\lambda(s)=\psi'_\lambda(\kappa(s))=\la\pi_\lambda^K(\kappa
(s)).v_\lambda,v_\lambda\ra\ \ \hbox{for all $s\in V^-V_0V^+$}.\leqno(2.3)$$

In (2.3) the second equality holds by definition so that it remains to
prove the first one. 
For each $Y\in\g_\C$ let  $L_Y$ the left invariant vector field on
$G_\C$ with $L_Y(\1)=Y$. The corresponding vectorfield on 
$V^-V_0V^+$ obtained by restriction and lifting is denoted by $\tilde
L_Y$. 
Note that the mapping 
$$\tilde L\:\g_\C\to\Der(C^\infty(V^-V_0V^+)),\ \ X\mapsto \tilde L_X$$
is a Lie algebra homomorphism and that all vectorfields 
$\tilde L_Y$ are holomorphic. 
Let $Y=Y_1+iY_2\in\p_\C^+$, $X\in U$ and $g:=\exp(X)\in G$. Since $v_\lambda$
is an analytic vector, we obtain that 
$$\eqalign{(\tilde L_Y.\psi'_\lambda)&(\Exp(X)) =
(\tilde L_{Y_1}.\psi'_\lambda)(\Exp(X))+i(\tilde L_{Y_2}.\psi'
_\lambda)(\Exp(X))\cr
&={d\over dt}\Big|_{t=0}\psi_\lambda(g\exp(tY_1))+i{d\over dt}\Big|_{t=0}
\psi_\lambda(g\exp(tY_2))\cr
&={d\over dt}\Big|_{t=0}\la\pi_\lambda(g\exp(tY_1)).v_\lambda,v_\lambda\ra
+i{d\over dt}\Big|_{t=0}\la\pi_\lambda(g\exp(tY_2)).v_\lambda,v_\lambda\ra\cr
&={d\over dt}\Big|_{t=0}\la\pi_\lambda(\exp(tY_1)).v_\lambda,
\pi_\lambda(g^{-1}).v_\lambda\ra
+i{d\over dt}\Big|_{t=0}\la\pi_\lambda(\exp(tY_2)).v_\lambda,
\pi_\lambda(g^{-1}).v_\lambda\ra\cr
&=\la d\pi_\lambda(Y_1).v_\lambda,v_\lambda\ra
+i\la d\pi_\lambda(Y_2).v_\lambda,v_\lambda\ra
=\la d\pi_\lambda(Y).v_\lambda,\pi_\lambda(g^{-1}).v_\lambda\ra=0.\cr}$$
Thus $(\tilde L_Y.\psi'_\lambda)\res_{\Exp(U)}=0$ and therefore 
$\tilde L_Y.\psi'_\lambda=0$ by the Identity Theorem for analytic
functions. Therefore $\psi'_\lambda$ is constant on all integral
curves of $\tilde L_Y$, i.e., 
$$\psi'_\lambda(p_-kp_+)=\psi'_\lambda(p_- k)\ \ 
\hbox{for all  $ p_-k p_+\in V^-V_0V^+$}.$$
Similarly one shows that $\psi'_\lambda$ is constant on the ``right
cosets'' of $V^-$, concluding the proof of (2.3). 

\par It follows from (2.3) that $\psi_\lambda$ extends to a
holomorphic function 
$$\tilde\psi_\lambda\:P^-\tilde{K_\C} P^+\to\C,\ \ s\mapsto
\la \pi_\lambda^K(\kappa(s)).v_\lambda,v_\lambda\ra.$$
In view of the Cartan-Weyl-Theorem of Highest Weight for finite dimensional 
representations of complex reductive Lie algebras, we have 
$HW(\Delta_k^+)=\z(\k_\C)^*+\Gamma$, where $\Gamma$ is an additive 
discrete subsemigroup of $i(\t\cap[\k,\k])^*$. Thus the continuity 
of $\tilde \Phi$ reduces to showing 
continuity of the maps
$$\tilde \Phi_\gamma\: \tilde{K_\C}\times \big(\z(\k_\C)^*+\gamma\big)\to \C, 
\ \ (k,\lambda)\mapsto \la \pi_\lambda^K(k).v_\lambda, v_\lambda\ra,$$ 
where $\gamma\in \Gamma$. 
Note that $\tilde {K_\C}\cong \z(\k_\C)\times [\tilde {K_\C},\tilde {K_\C}]
$ by the simple connectedness of $\tilde{K_\C}$.
Therefore an irreducible representation of $\tilde{K_\C}$ correspoding 
to $z+\gamma\in \z(\k_\C)^*+\gamma$  is a tensor product 
representation of a one-dimensional representation of $\z(\k_\C)$
with infinitesimal character $z$ and a highest weight representation 
of $[\tilde {K_\C},\tilde {K_\C}]$ corresponding to $\gamma$. 
This proves the continuity of the maps $\tilde \Phi_\gamma$ and 
completes the  proof of the proposition.  \qed

\Corollary II.21. The mappings $HW_\alpha\to (\hat S_\alpha,{\cal
T}_{hk}^\alpha)$ and 
$HW_\alpha\to (\hat S_\alpha,{\cal T}_G^\alpha)$ are continuous. 

\Proof. Recall that $\hat S_\alpha$ has a countable base (cf.\ Definition 
II.12). Thus by Theorem II.14 and [Dix82, Prop.\ 18.1.5],
it suffices to show  
$$\lambda_n\to\lambda\qquad  \hbox{in the euclidean topology }$$
implies 
$$\psi_{\lambda_n}\to\psi_\lambda\qquad  \hbox{compactly on $S$ and $G$}.$$
But this is immediate from Proposition II.20.\qed

Now we are going to prove that ${\cal T}_e^\alpha\subeq {\cal
T}_{hk}^e$. We start with a lemma 
which reduces the assertion to contraction representations.

\Lemma II.22. {\rm (Reduction to Contractions)} Let $\alpha\in {\cal
A}_h(S)$. Set $\g^\sharp \:=\g\oplus\R$ and $G^\sharp=G\times \R$. Then
the following assertions hold:
\item{(i)} The prescription 
$$W^\sharp\:=\{ (X,t)\in W\times\R^+\: \alpha(\Exp(iX))< e^t \}.$$
defines an open convex $\Ad(G^\sharp)$-invariant elliptic
cone in $\g^\sharp$ and a complex Ol'shanski\u\i{}
semigroup $S^\sharp\:=\Gamma_{G^\sharp}(W^\sharp)\subeq S\times \C$ can be build 
(cf.\ {\rm Definition II.2}).
\item{(ii)} If $(\pi,{\cal H})$ is an $\alpha$-bounded holomorphic representation
of $S$, then $\pi$ induces via
$\pi^\sharp(s,z)=e^{iz}\pi(s)$ 
a $\1$-bounded holomorphic representation of $S^\sharp$. Moreover, 
the prescription $\pi\mapsto \pi^\sharp$ defines a bijection between
$\hat S_\alpha$ and $\hat {S^\sharp}_\1$.
\item{(iii)} The $C^*$-algebras $C_h^*(S,\alpha)$ and $C_h^*(S^\sharp,\1)$ 
are isomorphic.

\Proof. (i) This is immediate from the definition of 
${\cal A}_h(S)$.
\par\nin (ii) {}From the construction of $S^\sharp$ and $\pi^\sharp$, 
it is clear that $\pi^\sharp$ is contractive, whenever $\pi$ is 
$\alpha$-bounded. Finally, recall from Lemma II.19(i) that every 
irreducible holomorphic representation of $S$, resp. $S^\sharp$, 
extends to a holomorphic representation of $\Gamma_G(W_{\rm max}^0)$,
resp. $\Gamma_{G^\sharp}(W_{\rm max}^0\oplus\R)$. In view of this, 
the second assertion is also clear.
\par\nin (iii) In the definition of  $C_h^*(S,\alpha)$ (cf.\
Definition I.3) we used the full cone  ${\cal P}_h(S,\alpha)$ of 
$\alpha$-bounded positive definite functions. However, in view of 
Segal's Theorem (cf.\ [Dix82, Lemma 2.10.1]), we may replace 
${\cal P}_h(S,\alpha)$  by the subcone of extremal generators 
$\Ext({\cal P}_h(S,\alpha))$. Now the assertion follows from (i), Lemma
II.6 and the construction of $C_h^*(S,\alpha)$ .\qed

A Lie group $G$ is called a {\it (CA)-Lie group} ({\it closed adjoint}) if 
$\Ad(G)$ is closed in $\Aut(\g)$. Note that all reductive and nilpotent 
Lie groups are (CA)-Lie groups.

\Lemma II.23.  Let $([\pi_{\lambda_n}])_{n\in\N}$ be a convergent 
sequence in $(\hat S_\alpha, {\cal T}_{hk}^\alpha)$ and 
$[\pi_{\lambda_0}]$ a limit point of it. Then the following 
assertions hold:

\item{(i)} The set $\{\lambda_n\: n\in\N\}$ is relatively compact in 
$i\t^*$. In particular, there exists a convergent subsequence 
$(\lambda_{n_k})$ with limit, say $\lambda_0'$. 

\item{(ii)} If, in addition, $G$ is a (CA)-Lie group, then we have 
$$\lambda_0'=\lambda_0+\mu_0$$
for some $\mu_0\in \N_0[\Delta^+]$, where $\N_0[\Delta^+]$ denotes
the additive subsemigroup of $i\t^*$ generated by $\Delta^+\cup\{0\}$.

\Proof. (i) W.l.o.g. we may assume that $S$ is simply connected. In view
of Lemma II.22, we may also assume that $\alpha\leq \1$.
\par If $[\pi_{\lambda_n}]\to [\pi_{\lambda_0}]$ in 
$(\hat S_{\alpha}, {\cal T}_{hk}^\alpha)$, then Theorem II.14, 
implies in particular that there exists unit vectors $v_n\in {\cal
H}_{\lambda_n}$, $n\in \N$, and a normalized 
highest weight vector $v_0$ of ${\cal H}_{\lambda_0}$ such that 

$$\la \pi_{\lambda_n}(s).v_n, v_n\ra \to \la \pi_{\lambda_0}(s).v_0,
v_0\ra
\qquad\hbox{compactly on $S$}.\leqno(2.4)$$

\par  For each 
$\lambda\in HW(G,\Delta^+)$ we write 
${\cal P}_\lambda\:={\cal  P}_{L(\lambda)}$ for the 
corresponding set of $\t_\C$-weights (cf.\ Definition II.18(a)). Note
that ${\cal P}_\lambda\subeq \lambda- \N_0[\Delta^+]$ since 
$(\pi_\lambda, {\cal H}_\lambda)$ is a highest weight 
representation with respect to $\Delta^+$ and highest weight
$\lambda$. 

\par Let $\t^+\:=\{X\in \t\: (\forall \alpha\in \Delta^+)\ i\alpha(X)>0\}$ and
let $X\in W\cap \t^+$. Then we have for all
$\lambda\in HW(G,\Delta^+)$ and $\mu_\lambda\in {\cal P}_\lambda$ 
that $\mu_\lambda(iX)\leq \lambda(iX)$. We show that 
$$\sup_{n\in \N_0} \lambda_n(iX) \leq 0\quad\hbox{and}\quad 
\inf_{n\in \N_0}\lambda_n(iX)>-\infty.\leqno(2.5)$$
The first assertion in (2.5) is immediate from our reduction to  
$\alpha\leq 1$. If the second assertion were false, we would find a
subsequence $(\lambda_{n_k})_{k\in \N}$ of 
$(\lambda_n)_{n\in\N}$ such that 
$\lim_{k\to\infty} \lambda_{n_k}(iX)=-\infty$. But this implies that 
$$0\leq \lim_{k\to\infty}\la \pi_{\lambda_{n_k}}(\Exp(iX)).v_{n_k},
v_{n_k}\ra \leq \lim_{k\to\infty} e^{\lambda_{n_k}(iX)} =0.$$
In view of (2.4), this means that $e^{\lambda_0(iX)}=
\la \pi_{\lambda_0}(\Exp(iX)).v_0, v_0\ra=0$;
a contradiction, completing the proof of (2.5).

\par By the definition of ${\cal A}_h(S)$ we have
$\alpha=\alpha_C$ for some closed convex subset $C\subeq \g^*$ with
$\oline {B(C)}=-\oline W$. In particular we have $X\in -\Int B(C)$.
Since $\lim C= -W^\star$, we see that $C$ is pointed (cf.\ Definition
II.4(c)), and so the evaluation mapping in $X$

$${\rm ev}_X\: C\to \R, \ \ \lambda\mapsto -\lambda(X)$$ is proper
(cf.\ [Ne99, Cor.\ V.1.16]). In view of Lemma II.19(ii), we have
$HW_\alpha\subeq -iC$, and so the properness of ${\rm ev}_X$ together
with (2.5) imply (i).

\par\nin (ii) For $\lambda\in HW(G,\Delta^+)$ and $v\in {\cal H}_\lambda$ we
write $v=\sum_{\mu\in {\cal P}_\lambda} v^\mu$ for the Fourier 
series of $v$ with respect to the $\t_\C$-action. 
For each $n\in \N$ and $\mu\in\N_0[\Delta^+]$ we set $a_{n, \mu}\:=\la
v_n^{\lambda_n -\mu}, v_n^{\lambda_n -\mu}\ra$. 
It follows from (2.4) that 
$$\la \pi_{\lambda_n}(\Exp(X)).v_n, v_n\ra \to \la
\pi_{\lambda_0}(\Exp(X)).v_0,v_0\ra
\quad\hbox{compactly on $\t+i(W\cap\t^+)$}.$$ 
Considering the
corresponding Fourier series, this means that 
$$e^{\lambda_n(X)}\sum_{\mu\in\N_0[\Delta^+]} a_{n, \mu}e^{-\mu(X)}
\to e^{\lambda_0(X)}$$
converges compactly on $\t+i(W\cap\t^+)$ or, equivalently, that 
$$\sum_{\mu\in \N_0[\Delta^+]} a_{n,\mu}e^{-\mu(X)}
\to e^{(\lambda_0-\lambda_0')(X)}\quad \hbox{compactly 
on $\t+i(W\cap\t^+)$}.\leqno(2.6)$$
\par Set $T\:=\exp(\t)$ and note that $\Ad(T)$ is a compact 
group since $G$ was supposed to be a (CA)-Lie group 
(cf.\ [Ne99, Sect.\ VII.1]). 
Thus (2.6) together 
with the Peter-Weyl Theorem applied to the action of $\Ad(T)$ on the 
Fr\'echet space $C^\infty(\t+i(W\cap\t^+))$ implies that
$a_{n,\mu}\to 0$ except for $\mu=\mu_0\:=\lambda_0'-\lambda_0$.
This completes the proof of (ii).\qed 

For $\lambda\in HW(G,\Delta^+)$ we also consider the $\k_\C$-module 
$F(\lambda)$  as $\k_\C+\p^+$-module  with trivial 
$\p^+$-action. Further 
we associate to $\lambda$ the {\it generalized Verma module}
$$N(\lambda)\:={\cal U}(\g_\C)\otimes_{{\cal U}(\k_\C+\p^+)} F(\lambda).$$ 
Note that $L(\lambda)$ is the unique irreducible 
quotient of the ${\cal U}(\g_\C)$-module $N(\lambda)$ 
(cf.\ [Ne99, Sect.\ IX.1]).

\Definition II.24. Let $\alpha\in {\cal A}_h(S)$ be an absolute 
value for $S$. Then we call $\alpha$ {\it generic} if $L(\lambda)=N(\lambda)$
holds for all $\lambda\in HW_\alpha$. \qed

\Proposition II.25. If  $G$ is a (CA)-Lie group and $\alpha\in 
{\cal A}_h(S)$ is generic, then the mapping 
$$(\hat S_\alpha, {\cal T}_{hk}^\alpha)\to HW_\alpha, \ \
 [\pi_\lambda]\mapsto \lambda$$
is continuous, i.e., ${\cal T}_e^\alpha\subeq {\cal T}_{hk}^\alpha$.

\Proof. We use the notaion of Lemma II.23 and its proof. 
Let $[\pi_{\lambda_n}]\to [\pi_{\lambda_0}]$ in 
$(\hat S_{\alpha}, {\cal T}_{hk}^\alpha)$. We have to show that 
$\lambda_n\to \lambda_0$. By Lemma II.23(i) we may assume that 
$(\lambda_n)_{n\in\N}$ converges to $\lambda_0'$. 
By Lemma II.23(ii) and its proof we find $\mu_0\in\N_0[\Delta^+]$
with $\lambda_0'=\lambda_o+\mu_0$, 
unit vectors $v_n\in L(\lambda_n)^{\lambda_n-\mu_0}$, $n\in\N$, 
and a normalized highest weight vector $v_0$ of $L(\lambda_0)$ such that 
(2.4) holds. 

\par  To complete the proof of the proposition, we have to show that 
$\mu_0=0$. Let  $Y\in \g_\C^\alpha$ for some $\alpha\in \Delta^+$. 
After multiplying $Y$ 
with a small non-zero scalar, we may assume that
there exists an open subset $U\subeq \t+ i(W\cap\t^+)$ such that 
$V=\Exp(-\oline Y)\Exp(U)\Exp(Y)\subeq P^-\tilde{K_\C}P^+\cap S$.
Now it follows from (2.4)that 
$$\eqalign{&\la \pi_{\lambda_n}(\Exp(-\oline Y)\Exp(X)\Exp(Y)).v_n, 
v_n\ra =\la \pi_{\lambda_n}(\Exp(X)\Exp(Y)).v_n, 
\pi_{\lambda_n}(\Exp(Y)).v_n\ra\cr
&\quad =\sum_{k=0}^\infty e^{(\lambda_n-\mu_0 +k\alpha)(X)}
{1\over (k!)^2} \la d\pi_{\lambda_n}(Y)^k.v_n, d\pi_{\lambda_n}(Y)^k.v_n\ra\cr
&\quad \to\la\pi_{\lambda_0}(\Exp(-\oline Y)\Exp(X)\Exp(Y)).v_0, v_0\ra
=e^{\lambda_0(X)}\cr}\leqno(2.7)$$ 
converges compactly on $U$. Again by the Peter-Weyl Theorem (cf.\ the proof
of Lemma II.23(ii)) we deduce  
that 
$$(\forall Y\in \n\:=\bigoplus_{\alpha\in\Delta^+}\g_\C^\alpha)\quad 
\lim_{n\to\infty} \|d\pi_{\lambda_n}(Y).v_n\|=0.\leqno(2.8)$$ 

\par As $\lambda_n\to \lambda_0'$ we may assume that 
$\lambda_n\res_{\t_\C\cap [\k_\C,\k_\C]}=
\lambda_0'\res_{\t_\C\cap [\k_\C,\k_\C]}$ for all $n\in \N$. Set 
$\Lambda\:=\{\lambda_n\: n\in \N\}\cup\{\lambda_0'\}$ and note that 
$\Lambda\subeq HW_\alpha$ since $HW(G,\Delta^+)$ is closed in $i\t^*$
(cf.\ [Kr99a, Sect.\ IV]) and 
$HW_\alpha$ is closed in $HW(G,\Delta^+)$ (cf.\ Lemma II.19(ii)). 

\par As $\alpha$ is generic, [Kr99b] implies that we can 
identify all $L(\lambda)=N(\lambda)$, 
$\lambda\in \Lambda$, with $L(\lambda_1)$ as $[\k_\C,\k_\C]+\p^+$-modules.
Within this identification all operators 
$d\pi_\lambda(Y)$, $\lambda\in\Lambda$, 
coincide for $Y\in \p^+$. Further the scalar products of the various
$L(\lambda)$, $\lambda\in\Lambda$, define a family of inner 
products $(\la\cdot,\cdot\ra_\lambda)_{\lambda\in \Lambda}$ on 
$L(\lambda_1)$. Also from [Kr99b] we can deduce that for all
$v\in L(\lambda_1)$ the mapping 
$$\Lambda\to]0,\infty[, \ \ \lambda\mapsto\la v,v\ra_\lambda$$
is continuous. In particular we may assume that $(v_n)_{n\in\N}$ converges in 
$L(\lambda_1)^{\lambda_1-\mu_0}$ with limit $v_0'\neq 0$. 
Thus (2.8) gives for all $Y\in\n$ 
$$\|d\pi_{\lambda_0'}(Y).v_0'\|=
\|d\pi_{\lambda_1}(Y).v_0'\|_{\lambda_0'}=
\lim_{n\to\infty}\|d\pi_{\lambda_1}(Y).v_n\|_{\lambda_n}=
\lim_{n\to\infty}\|d\pi_{\lambda_n}(Y).v_n\|=0,$$
where the subscripts at the various norms indicate that we have identified 
the corresponding $L(\lambda)$ with $L(\lambda_1)$. 
Thus $v_0\in L(\lambda_0')^{\lambda_0'-\mu_0}$ is a primitive 
element and so $\mu_0=0$ as was to be shown. \qed

We now give an example that Proposition II.25 becomes false if 
$\alpha$ is non generic.

\Example II.26. Let  $G:=\tilde\Sl(2,\R)$
be the universal covering group of $\Sl(2,\R)$. We choose
$$U=\pmatrix{0 & 1\cr -1 & 0\cr},\qquad T=\pmatrix{0 & 1\cr 1 & 0\cr},\quad\hbox{and}\quad
H=\pmatrix{1 & 0\cr 0 & -1\cr} $$
as a basis for $\g:=\sL(2,\R)$.
Then  $\t:=\R U$ is a compactly embedded Cartan  subalgebra. Let $\alpha\in
i \t^*$ be defined by $\alpha(U)=-2i$. The root system of $\g$ is given by 
$\Delta=\{ \pm\alpha\}$ with root spaces 
$\g_\C^\alpha=\C(T+iH)$ and $\g_\C^{-\alpha}
=\C(T-iH)$. We define a positive system by $\Delta^+:=\{\alpha\}$. We
denote by $\kappa$ the Cartan-Killing form of $\g$. 
Then the upper light cone 
$$W :=\{X=uU+tT+hH:u\geq 0,\ \kappa(X,X)\leq 0\}
=\{X=uU+tT+hH:u\geq 0,\ h^2+t^2-u^2\leq 0\}$$
is an invariant pointed cone in $\g$. Moreover, $W$ is up to sign the unique
invariant elliptic cone in $\g$ (cf.\ [HiNe93, Th.\ 7.25]). Let
$S:=\Gamma_G(W)$ be the complex Ol'shanski\u\i{}  semigroup 
corresponding to $G$ and $W$. 

\par In the following we identify $\t_\C^*$ with $\C$ via the isomorphism 
$\t_\C^*\to \C, \ \lambda\mapsto \lambda(iU)$.  
Then $HW(G,\Delta^+)=]-\infty, 0]$. 

\par If $\alpha\in {\cal A}_h(S)$ is an absolute value, then Lemma
II.19(ii) implies that $HW_\alpha=]-\infty, t_\alpha]$ for some
$t_\alpha\leq 0$. Conversely it follows from [Ne99, Th.\ VIII.3.21] that for
each $t_\alpha\leq 0$ there exists an $\alpha\in {\cal A}_h(S)$  with 
$HW_\alpha=]-\infty, t_\alpha]$. Thus 
${\cal A}_h(S)\to ]-\infty, 0], \ \alpha\mapsto t_\alpha$ is bijective, hence
gives a parametrization of ${\cal A}_h(S)$ with non-positive real
numbers. The generic absolute values correspond to $]-\infty,0[$ and the only 
non-generic one is $\alpha=\1$ which corresponds to $t_\alpha=0$ (cf.\ [Kr98, Ex.\ III.7]).

\par Let us now show that Proposition II.25 becomes false for 
the non-generic absolute value $\alpha=\1$. 
We choose $\lambda_n=-{1\over n}$, $n\in\N$, and show that 
$[\pi_{-{1\over n}}]\to [\pi_{-1}]$ in 
$(\hat S_\1, {\cal T}_{hk}^\1)$. Since we evidently have 
$[\pi_{-{1\over n}}]\to [\pi_0]$, this will give us the non-continuity of the 
map in Proposition II.25.

\par For each $n\in \N$ Let $v_n\in 
L(-{1\over n})^{-{1\over n}-2}=L(\lambda_n)^{\lambda_n-\alpha}$, 
$n\in\N$, be a normalized vector. Further let $v_{-1}\in L(-1)$ be 
a normalized highest weight vector. Then by Theorem II.14 we will have $[\pi_{-{1\over n}}]\to [\pi_{-1}]$ in 
$(\hat S_\1, {\cal T}_{hk}^\1)$ if $\la\pi_{-{1\over n}}(s).v_n, v_n\ra\to \la\pi_{-1}(s).v_{-1}, v_{-1}\ra$
holds uniformly on compact subsets $U\subeq S$. 

\par Set $E\:=T+iH\in \g_\C^\alpha$. Now for $U\subeq S$ compact we 
find an $R>0$ such that 
$U\subeq V^- W V^+$ with $W\subeq \tilde{K_\C}$ compact,  
$V^-=\{\Exp(u\oline E)\: u\in \C, \ |u|\leq R\}$ and    
$V^+=\{\Exp(wE)\: w\in \C,\ |w|\leq R\}$. Further each $s\in S$ can uniquely be written as 
$s=\Exp(-u(s)\oline E)\kappa(s) \Exp(w(s) E)$. Thus we get that 
$$\eqalign{\la\pi_{-{1\over n}}(s).v_n, v_n\ra&=
\la\pi_{-{1\over n}}(\Exp(-u(s)\oline E)\kappa(s)\Exp(w(s)E)).v_n, v_n\ra\cr 
&=\kappa(s)^{\lambda_n-\alpha} + w(s)\oline {u(s)}
\la d\pi_{-{1\over n}}(E).v_n, d\pi_{-{1\over n}}(E).v_n\ra.\cr}\leqno(2.9)$$
Now the formulas 
in [La85, Ch.\ VI] show that 
$$\la d\pi_{-{1\over n}}(E).v_n, d\pi_{-{1\over n}}(E).v_n\ra=
{\big({1\over n}\big)^2\over \big( 2+{1\over n}\big)}\to 0.\leqno(2.10)$$
Thus we conclude from (2.9) and (2.10) that 
$$\la\pi_{-{1\over n}}(s).v_n, v_n\ra\to \kappa(s)^{-\alpha}=\la\pi_{-1}(s).v_{-1}, v_{-1}\ra$$
uniformly on $U$ as was to be shown. \qed

\subheadline{The Borel structure on the dual}

In this subsection we investigate the Borel structures on 
$\hat S_\alpha$ induced from our various topolgies on the dual.

\Lemma II.27. Let  ${\cal A}$ be a $C^*$-algebra and 
$x\in{\cal A}_+:=\{y\in{\cal A}\: y\geq 0\}$. 
Then the mapping 
$$\hat {\cal A}\to\R,\ \ [\pi]\mapsto\sup\Spec(\pi(x))=\|\pi(x)\|$$
lower semicontinuous. 

\Proof. [Dix82, Prop.\ 3.3.2]. \qed

\Proposition II.28. Let ${\cal B}_{hk}^\alpha$ be the Borel structure on 
$\hat S_\alpha$ induced from ${\cal T}_{hk}^\alpha$ and similarily 
${\cal B}_e^\alpha$ the one induced from ${\cal T}_e^\alpha$.
Then the mapping $(\hat S_\alpha,{\cal B}_{hk}^\alpha)\to (\hat S_\alpha,
{\cal B}_e^\alpha)$ is measurable, i.e., ${\cal B}_e^\alpha\subeq
{\cal B}_{hk}^\alpha$.

\Proof. Let  $X\in W\cap \t^+$. Then $\Exp(iX)$ is a symmetric element 
of $S$ and so $j(\Exp(iX))=j(\Exp(i{1\over 2}X))^2$ is  positve in 
$C_h^*(S,\alpha)$.
Let $v_\lambda\in{\cal H}_\lambda$
denote a  normalized highest weight vector. Since 
${\cal P}_\lambda\subset \lambda-\N_0[\Delta^+]$ we get 
$$\eqalign{\sup\Spec(\tilde{\pi_\lambda}(j(\Exp(iX)))&=
\sup\Spec(\pi_\lambda(\Exp(iX))) \cr &=\la\pi_\lambda(\Exp(iX)).v_\lambda,
v_\lambda\ra=e^{i\lambda(X)}.\cr}$$
In view of Lemma II.27,  for all $\beta\in \R$ 
the subsets 
$$\eqalign{I_{X,\beta}:=&\{\lambda\in\hat S_\alpha: 
e^{i\lambda(X)}\leq e^\beta\}\cr =&\{\lambda\in\hat S_\alpha:i\lambda(X)\leq
\beta\}\cr}$$
are closed in $(\hat S_\alpha, {\cal T}_{hk}^\alpha)$. 
Now the assertion of the proposition follows since the 
system $\{I_{X,\beta}\: X\in W\cap\t^+,\ \beta\in\R\}$ generates 
the euclidean Borel structure on $\hat S_\alpha$.\qed

\subheadline{The main results}

\Theorem II.29. {\rm (The topologies on $\hat S_\alpha$)}
Let $S=\Gamma_G(W)$ be a complex Ol'shanski\u\i{}
semigroup, $\alpha\in {\cal A}_h(S)$ an absolute value on $S$  and  
$\hat S_\alpha$ the equivalence classes of all $\alpha$-bounded irreducible
holomorphic representations of $S$. Let ${\cal T}_{hk}^\alpha$ denote
the hull-kernel topology on $\hat S_\alpha$ obtained from
$C_h^*(S,\alpha)$, further ${\cal T}_G^\alpha$ the topology induced from
$\hat G$ and ${\cal T}_e^\alpha$ the euclidean topology obtained from
the parametrization with heighest weights, and write 
${\cal B}_{hk}^\alpha$, ${\cal B}_G^\alpha$ and ${\cal B}_e^\alpha$
for the corresponding Borel structures. 

\item{(i)} We have 
$${\cal T}_{hk}^\alpha\subeq {\cal T}_G^\alpha\subeq{\cal T}_e^\alpha\quad 
\hbox{and}\quad {\cal B}_{hk}^\alpha={\cal B}_G^\alpha={\cal B}_e^\alpha.$$

\item{(ii)} If, in addition, $G$ is a (CA)-Lie group and 
$\alpha$ is generic, then 
$${\cal T}_{hk}^\alpha={\cal T}_G^\alpha={\cal T}_e^\alpha.$$

\Proof. (i) The inclusion for the topologies follows from 
Proposition II.16 and Corollary II.21. 
Finally, the identity for the Borel structures follows from the 
chain of inclusions for the topologies together with 
Proposition II.28.
\par\nin (ii) In view of (i), this follows from Proposition 
II.25.\qed

Now we give two applications of Theorem II.24 to the structure of
$C_h^*(S,\alpha)$  and the abstract representation theory of complex 
Ol'shanski\u\i{} semigroups. In the following two statements we use
the language of [Dix82].

\Corollary II.30. Let $G$ be a (CA)-Lie group, 
$\alpha$ a generic  absolute value of $S$ and ${\frak A}_\alpha$ the
$C^*$-algebra defined by the continuous field $\big({\cal K}({\cal
H}_\lambda)\big)_{\lambda\in HW_\alpha}$ of elementary $C^*$-algebras. Then the
mapping 
$$C_h^*(S,\alpha)\to {\frak A}_\alpha, \ \ x\mapsto (\tilde
\pi_\lambda(x))_{\lambda
\in HW_\alpha}$$
(cf.\ {\rm Lemma II.6}) is an isomorphism of $C^*$-algebras. 

\Proof. In view of Theorem II.8 and Theorem II.29(ii), $C_h^*(S,\alpha)$ is
a CCR $C^*$-algebra with Hausdorff spectrum, and so the assertion
follows from [Dix82, Th.\ 10.5.4].\qed

Recall the notion of multiplicity free representations: A holomorphic
representation $(\pi, {\cal H})$ of $S$ is called {\it multiplicity
free} if its commutant $\pi(S)'$ in $B({\cal H})$ is abelian. 
Even though a more general formulation of the next result is possible, we
emphasize on multiplicity free representations, since in the author's 
opinion the most interesting 
examples of holomorphic representations of complex Ol'shanski\u\i{} 
semigroups are multiplicity free.

\Corollary II.31. Let $(\pi, {\cal H})$  be a holomorphic 
multiplicity free representation of $S$ and $\alpha$ be the absolute 
value defined by $\alpha(s)=\|\pi(s)\|$. Then there exists a Radon
measure $\mu$ on the euclidean space $HW_\alpha\subeq i\t^*$  and a direct
integral of representations
$$\Big(\int_{HW_\alpha}^\oplus\pi_\lambda\
d\mu(\lambda),\int_{HW_\alpha}^\oplus {\cal H}_\lambda\
d\mu(\lambda)\Big)$$
such that $(\pi, {\cal H})$ is unitarily 
equivalent with $\big(\int_{HW_\alpha}^\oplus\pi_\lambda\
d\mu(\lambda),\int_{HW_\alpha}^\oplus {\cal H}_\lambda\
d\mu(\lambda)\big)$.

\Proof. Recall that the $\alpha$-bounded holomorphic representations 
can be modelled by the representations of the CCR $C^*$-algebra $C_h^*(S,\alpha)$
(cf.\ Lemma II.6). In view of this, the assertion follows from 
[Dix82, Th.\ 8.6.5] and Theorem II.29(i).  \qed

\subheadline{Continuous traces}

In this final subsection we give a sufficient criterion for 
$C_h^*(S,\alpha)$ to have  continuous trace. We will illustrate this
criterion in various examples. 

Recall that for each holomorphic highest weight representation 
$(\pi_\lambda, {\cal H}_\lambda)$ of $S$ all operators
$\pi_\lambda(s)$, $s\in S$, are of trace class (cf.\ [Ne99, Th.\
XI.6.1]), and so the notion

$$\Theta_\lambda\: S\to \C, \ \ s\mapsto \tr\pi_\lambda(s)$$
makes sense. We call $\Theta_\lambda$ the {\it character} of
$(\pi_\lambda, {\cal H}_\lambda)$ and note that $\Theta_\lambda$ is a
holomorphic function on $S$ (cf.\ [Ne99, Prop.\ XI.6.4]).

\Definition II.32. (cf.\ [Dix82, Sect.\ 4.5.2]) Let ${\cal A}$ be a $C^*$-algebra and ${\cal A}_+$
the cone of positive elements in ${\cal A}$. Let $\p\subeq {\cal A}_+$ be the subcone 
of those elements of $x\in {\cal A}_+$ for which the mapping 
$\hat {\cal A}\to[0,\infty], \ [\pi]\mapsto \tr\pi(x)$ is 
finite and continuous. Then $\p$ is the positive portion of a
two-sided ideal $\m$ of ${\cal A}$. We say that ${\cal A}$ has
{\it continuous trace} if $\m$ is dense in ${\cal A}$. \qed

\Proposition II.33. Let $S$ be a complex Ol'shanski\u\i{} semigroup 
and $\alpha\in{\cal A}_h(S)$ be an absolute value on it. 
Then the following assertions hold:

\item{(i)} If there exists an open subset $U\subeq W\cap \t$
such that for each $X\in U$ the mapping 
$$\phi_X\: HW_\alpha\to \R^+, \ \ \lambda\mapsto
\Theta_\lambda(\Exp(iX))$$
is continuous, then $C_h^*(S,\alpha)$ has continuous trace. 

\item{(ii)} If $\alpha$ is generic, then $C_h^*(S,\alpha)$ has continuous
trace.

\Proof. (i)  Let $V\:=\Ad(G).U$ and note that $V$ is an open subset of
$W$. Let $\p$ and $\m$ as in Definition II.29. Then $\p\supeq
j(\Exp(iU))$  by Theorem II.24 and so $\p\supeq j(\Exp(iV))$. Since every holomorphic
representation of $S$ which vanishes on $\Exp(iV)$ has to be
constant, it follows that $j(\Exp(iV))$ generates a dense ideal  of
$C_h^*(S,\alpha)$ (cf.\ Lemma II.6). In particular $\m$ is dense as was
to be shown. 
\par\nin (ii) In view of (i), this follows from [Kr99a, Lemma IV.9].\qed

\Example II.34. (a) Suppose that $G$ is a simply connected hermitian
Lie group. Set $\t_0=\t\cap [\k, \k]$ and note that $\t=\z(\k)\oplus
\t_0$. According to this decomposition, every element $\lambda\in
i\t^*$ can be written as $\lambda=\lambda_z+\lambda_0$ with $\lambda_z\in 
i\z(\k)^*$ and $\lambda_0\in i\t_0^*$. 
\par Assume that the absolute value $\alpha$ satisfies the following
condition:

\item{(CT)} For each $\lambda\in HW_\alpha$ which is not isolated there exits
an $\eps>0$ such that $]1-\eps, 1+\eps[\lambda_z+\lambda_0 \subeq HW(G,\Delta^+)$. 

\par\nin Note that this condition excludes the first reduction
points in $HW(\tilde G, \Delta^+)$ and (CT) is exactly the condition 
for $\alpha$ being generic  (cf.\ [Ne99, Ch.\ X]). 
Therefore Proposition II.33(ii) applies and shows that $C_h^*(S,\alpha)$ has
continuous trace. 

\par\nin (b) We now discuss (a) in the special case of $G:=\tilde\Sl(2,\R)$
the universal covering group of $\Sl(2,\R)$. We use the notation of 
Example II.27. 

\par Recall that $\t_\C^*$ was identified with  $\C$ and within this 
identification we have $HW(G,\Delta^+)=]-\infty, 0]$. 
Then for all $u\in\R^+$ one has 

$$\Theta_\lambda(\Exp(iuU))=\cases{{e^{u\lambda}\over
1-e^{2u}} & for $\lambda<0$\cr 
1 & for $\lambda=0$\cr}\leqno(2.11)$$ 
(this is a special case of [Kr99a, Lemma IV.9]; see also [La85, Ch.\
VII,\S 4, Th.\ 5]).
\par Also recall our identifaction for the absolute values 
${\cal A}_h(S)\to ]-\infty, 0], \ \alpha\to t_\alpha$. 
Now it follows from (a) and (2.11) that $C_h^*(S,\alpha)$ has continuous
trace if and only if $t_\alpha<0$, i.e., $\alpha$ is generic. \qed

\def\entries{

\[Dix82 Dixmier, J., ``$C^*$-{\it Algebras},'' North Holland, Amsterdam,
New York, Oxford, 1982

\[HiNe93 Hilgert, J. and K.--H. Neeb,  ``Lie semigroups and their Applications,'' 
Lecture Notes in Math.\ {\bf 1552}, Springer, 1993

\[Kr98  Kr\"otz, B.,  {\it On Hardy and Bergman spaces on complex Ol'shanski\u\i{}
semigroups}, Math.\ Ann. {\bf 312} (1998), 13--52 

\[Kr99a ---, {\it The Plancherel Theorem for Biinvariant Hilbert Spaces}, Publ. RIMS, Kyoto University, 
{\bf 35(1)} (1999), 91--122, 

\[Kr99b ---, {\it Norm inequalities for unitarizable highest weight modules}, 
Ann. Inst. Four. {\bf 49(4)} (1999), 1241--1264

\[La85 Lang, S., ``SL(2)", GTM {\bf 105}, Springer, 1985

\[Ne99 Neeb, K.--H, ``Holomorphy and Convexity in Lie Theory,'' 
Expositions in Mathematics, de Gruyter, in press

\[SW75 Stein, E., and Weiss, G., ``Introduction to Fourier Analysis on Euclidean
Spa\-ces,'' Princeton University Press, Princeton (NJ), 1975}

{\sectionheadline{\bf References}
\frenchspacing
\entries\par}
\dlastpage
\vfill\eject
\bye